\documentclass[12pt,reqno]{amsart}
\usepackage[headings]{fullpage}
\usepackage{amsmath,amssymb,amsthm}
\usepackage[bookmarks=true,%
colorlinks=true,%
linkcolor=blue,%
citecolor=blue,%
filecolor=blue,%
menucolor=blue,%
urlcolor=blue,%
breaklinks=true]{hyperref}
\usepackage{graphicx,svg,subcaption,url}
\usepackage{fullpage,comment}
\usepackage[shortlabels]{enumitem}
\usepackage{mathrsfs,mathtools,scalerel}
\usepackage{tikz,tikz-cd}
\usetikzlibrary{knots,math,decorations.markings,decorations.pathreplacing,patterns}
\usepackage{ifthen}

\allowdisplaybreaks

\makeatletter
\renewcommand*{\fps@figure}{htpb!}
\makeatother

\newlist{enumcc}{enumerate}{1}
\setlist[enumcc]{leftmargin=0pt, itemindent=*}
\setlist{font=\normalfont}

\newtheorem{theorem}{Theorem}[section]
\theoremstyle{definition}
\newtheorem{proposition}[theorem]{Proposition}
\newtheorem{lemma}[theorem]{Lemma}
\newtheorem{definition}[theorem]{Definition}
\newtheorem{remark}[theorem]{Remark}
\newtheorem{corollary}[theorem]{Corollary}

\newtheorem{question}[theorem]{Question}

\newtheorem{example}[theorem]{Example}

\tikzset{knot diagram/every knot diagram/.style={
    background color=gray!20,clip width=7,end tolerance=5pt,clip radius=0.1cm}}
\newenvironment{linkdiag}[1][0.5]{\mathop{}\!
  \begin{tikzpicture}[scale=0.8,baseline=(ref.base)]
    \node (ref) at (0.5,{#1}){\phantom{$-$}};}{\end{tikzpicture}\!\mathop{}}
\tikzset{-o-/.code 2 args={
  \pgfkeysalso{decoration={markings,mark=at position #1 with {\arrow{#2}}},
    postaction={decorate}}
}}

\providecommand{\abs}[2][]{#1\lvert{#2}#1\rvert}
\newcommand{\bracket}[2][]{#1\langle{#2}#1\rangle}

\DeclareMathOperator{\id}{id}
\DeclareMathOperator{\Span}{span}

\DeclareMathOperator{\Hom}{Hom}
\DeclareMathOperator{\End}{End}

\DeclareMathOperator{\Inv}{Inv}
\DeclareMathOperator{\ev}{ev}
\DeclareMathOperator{\coev}{coev}
\DeclareMathOperator{\RT}{RT}

\newcommand{\nats}{\mathbb{Z}_{\ge0}}
\newcommand{\ints}{\mathbb{Z}}
\newcommand{\rats}{\mathbb{Q}}

\newcommand{\groundR}{R_\colorL}
\newcommand{\groundK}{k_\colorL}

\newcommand{\CC}[1]{\{#1\}}

\newcommand{\embed}{\hookrightarrow}
\newcommand{\onto}{\twoheadrightarrow}

\newcommand{\gsl}[1][2|1]{{\mathfrak{sl}_{#1}}}
\newcommand{\Uq}{\mathrm{U}_q}
\newcommand{\Uqsl}[1][2|1]{\Uq(\gsl[#1])}
\newcommand{\Repsl}[1][2|1]{\Rep_q(\gsl[#1])}

\def\BZ{\mathbb Z}

\def\BQ{\mathbb Q}

\def\be{\begin{equation}}
\def\ee{\end{equation}}
\def\z{\zeta}
\def\fsl21{\gsl}
\def\Sk{\skein}

\newcommand{\LG}{\mathrm{LG}}
\newcommand{\GG}{\mathrm{\Gamma G}}
\newcommand{\skein}{\mathrm{Sk}}
\newcommand{\colorL}{\Lambda}
\newcommand{\Rep}{\mathsf{Rep}}
\newcommand{\cat}{\mathsf{C}}
\newcommand{\webcat}{{\mathsf{Web}}}
\newcommand{\ribcat}{\mathsf{Rib}}
\newcommand{\skeinballcat}{{\mathsf{Sk}(B^3)}}
\newcommand{\trees}{\mathfrak{T}}
\newcommand{\paths}{\mathfrak{P}}

\newcommand{\stsize}{\tiny}
\newcommand{\webnode}[4][9pt]{
  \node[circle,draw,fill=white,minimum size=#1,inner sep=0pt] (#3) at (#2) {};
  \node at (#3) {\tiny$#4$} % no ;
}
\newcommand{\Iskein}[2]{
\begin{linkdiag}
\fill[gray!20] (-0.3,-0.3)rectangle(1.3,1.3);
\webnode{0.5,0.1}{A}{#1};
\webnode{0.5,0.9}{B}{#2};
\draw[thick] (A) -- (-0.2,-0.3) (A) -- (1.2,-0.3)
  (B) -- (-0.2,1.3) (B) -- (1.2,1.3)
  (A) -- (B);
\end{linkdiag}
}
\newcommand{\Hskein}[2]{
\begin{linkdiag}
\fill[gray!20] (-0.3,-0.3)rectangle(1.3,1.3);
\webnode{0.1,0.5}{A}{#1};
\webnode{0.9,0.5}{B}{#2};
\draw[thick] (A) -- (-0.2,-0.3) (B) -- (1.2,-0.3)
  (A) -- (-0.2,1.3) (B) -- (1.2,1.3)
  (A) -- (B);
\end{linkdiag}
}
\newcommand{\iddiag}{
\begin{linkdiag}
\fill[gray!20] (-0.3,-0.3)rectangle(1.3,1.3);
\begin{scope}[thick]
\draw (-0.2,1.3)..controls(0.2,0.5)..(-0.2,-0.3);
\draw (1.2,1.3)..controls(0.8,0.5)..(1.2,-0.3);
\end{scope}
\end{linkdiag}
}

\begin{document}
\title[Skein theory for $\mathfrak{sl}_{2|1}$-webs]{
  Skein theory for $\mathfrak{sl}_{2|1}$-webs}

\author{Stavros Garoufalidis}
\address{% Max Planck Institute for Mathematics \\
% Vivatsgasse 7, 53111 Bonn, GERMANY \newline
Shenzhen International Center for Mathematics, Department of Mathematics \\
Southern University of Science and Technology \\
1088 Xueyuan Avenue, Shenzhen, Guangdong, China \newline
{\tt \url{http://people.mpim-bonn.mpg.de/stavros}}}
\email{stavros@mpim-bonn.mpg.de}

\author{Tao Yu}
\address{Shenzhen International Center for Mathematics \\
Southern University of Science and Technology \\
1088 Xueyuan Avenue, Shenzhen, Guangdong, China}
\email{yut6@sustech.edu.cn}

\thanks{
  {\em Keywords and phrases}:
  skein theory, 3-manifolds, Lie superalgebras, $\mathfrak{sl}_{2|1}$, 
  webs, cabling, Hecke algebras, braids, Links--Gould polynomial, 
  multivariable Alexander polynomial.
}

\date{27 November 2026}%{}\today}
% \dedicatory{}

\begin{abstract}
  We develop an equivalence of the skein category of 
  representations of $\Uqsl$ with a diagrammatic category of 
  $\mathfrak{sl}_{2|1}$-webs, a dimension count, a basis and  
  a reduction algorithm for the skein algebra of the 3-ball. 
  We study the skein algebra of the annulus and prove cabling
  formulas for the colored Links--Gould polynomial of links,
  and we identify the $q=1$ specialization of the colored Links--Gould
  polynomial with the square of the Conway potential function.  
\end{abstract}

\maketitle

{\footnotesize
\tableofcontents
}

%%%%%%%%%%%%%%%%%%%%%%%%%%%%%%%%%%%%%%%%%%%%%%%%%%%%%%%%%%%%%%%%%%%%%%%%%%%%
%%%%%%%%%%%%%%%%%%%%%%%%%%%%%%%%%%%%%%%%%%%%%%%%%%%%%%%%%%%%%%%%%%%%%%%%%%%%

\section{Introduction}
\label{sec.intro}

In their pioneering work twenty years ago~\cite{GP:multi}, Geer--Patureau
introduced multivariable polynomials of links and knotted trivalent graphs 
using
the representation theory of $\gsl$. The feature of these invariants, new
at the time, was the non-semisimple aspect of the representation category 
which
lead to invariants that depended on continuous variables, other than the 
usual
variable $q$ which is omnipresent in quantum topology and the discrete 
variables
that correspond to weights of irreducible representations of simple Lie
algebras. This theory was later extended to a TQFT~\cite{BCGP}.

Our aim is to present a diagrammatic equivalence of this skein theory that
builds on, and unifies, the foundational work of
Links--Gould~\cite{Links-Gould}, Geer--Patureau~\cite{GP:multi},
Marin--Wagner~\cite{Marin-Wagner} and Anghel~\cite{Anghel}. It may be used 
as a
step towards the categorification of the $\fsl21$ knot polynomials, something
that we will discuss in subsequent work, as well as a way to define a
Habiro-cohomology valued TQFT, which we plan to discuss in a later 
publication. 

Turaev~\cite{Tu:conway} and Przytycki~\cite{Przytycki} introduced skein theory
for $\gsl[2]$ in the early 1990s, in the spirit of Jones and
Kauffman~\cite{Jones,Kauffman}. The famous Kauffman bracket relations can be
explained by the fact that $\Inv_{\Uqsl[2]}(V^{\otimes 4})$ is a 2-dimensional
space, where $V$ is the 2-dimensional fundamental representation of $\gsl[2]$.

Kuperberg~\cite{Kuperberg} extended skein theory to rank 2 Lie algebras. In 
the
case of $\mathfrak{sl}_3$, the corresponding diagrams (called webs) were
trivalent since the invariant space $\Inv_{\Uqsl[3]}(V^{\otimes 3})$ is
1-dimensional, where $V$ is the 3-dimensional fundamental representation of
$\mathfrak{sl}_3$.

In our case, the trivalent condition of our $\gsl$-webs comes from the fact 
that
the invariant space $\Inv_{\Uqsl}\big(\bigotimes_{i=1}^3 V(0,a_i)\big)$ is at
most 1-dimensional, where $V(0,a),a\ne0,-1$ is a family of 4-dimensional
irreducible representations indexed by a variable $a$; see
Equation~\eqref{Vten}.

\subsection{$\gsl$-webs}

To begin with, we choose a color lattice $\colorL$, a free abelian group of
finite rank, with a direct summand identified with $\ints$. In particular,
non-integral rational numbers are not in $\colorL$ by assumption.

Let $M$ be an oriented 3-manifold, and $V\subset\partial M$ be a finite set of
framed and typically colored points. Here, the framing is a choice of a nonzero
vector tangent to $\partial M$ at each point, and the color is an element of
$\colorL\setminus\{0,-1\}$. (These elements are called typical colors.)

A trivalent web in $(M,V)$ is a smoothly embedded graph $T$ whose vertices 
have
degree 1 or 3 such that the 1-valent vertices agree with $V$. By convention, 
the
empty set is \emph{not} a web even if $V=\emptyset$. A 
$\gsl$-web in $(M,V)$ is
a framed and typically colored trivalent web. The framing is a normal vector
field agreeing with that of $V$, which determines a cyclic order of edges 
around
each 3-valent vertex. The color is a choice of an orientation and a typical
element $a\in\colorL\setminus\{0,-1\}$ on each edge subject to the 
equivalence
\begin{equation}\label{eq-color-flip}
\begin{linkdiag}
\fill[gray!20] (0,0) rectangle (1.5,1);
\draw[thick,-o-={0.3}{<}] (0,0.4) -- node[pos=0.4,above]{$a$} (1.5,0.4);
\end{linkdiag}
=
\begin{linkdiag}
\fill[gray!20] (0,0) rectangle (1.5,1);
\draw[thick,-o-={0.75}{>}] (0,0.4) -- node[pos=0.6,above]{$a^\ast$} 
(1.5,0.4);
\end{linkdiag},\qquad
a^\ast:=-a-1,
\end{equation}
two possible conservation rules at each vertex, indicated by a sign $\pm$,
\begin{equation}\label{eq-vertex-sum}
a+b+c=-2\text{ in }
\begin{linkdiag}
\fill[gray!20] (-0.3,-0.3)rectangle(1.3,1.3);
\webnode{0.5,0.5}{V}{+};
\begin{scope}[thick]
\draw[-o-={0.8}{>}] (V) --node[pos=0.2,left]{$a$} (-0.2,1.3);
\draw[-o-={0.8}{>}] (V) --node[pos=0.2,right]{$c$} (1.2,1.3);
\draw[-o-={0.7}{>}] (V) --node[left]{$b$} (0.5,-0.3);
\end{scope}
\end{linkdiag},\quad\text{or}\quad
a+b+c=-1\text{ in }
\begin{linkdiag}
\fill[gray!20] (-0.3,-0.3)rectangle(1.3,1.3);
\webnode{0.5,0.5}{V}{-};
\begin{scope}[thick]
\draw[-o-={0.8}{>}] (V) --node[pos=0.2,left]{$a$} (-0.2,1.3);
\draw[-o-={0.8}{>}] (V) --node[pos=0.2,right]{$c$} (1.2,1.3);
\draw[-o-={0.7}{>}] (V) --node[left]{$b$} (0.5,-0.3);
\end{scope}
\end{linkdiag},
\end{equation}
and compatible with $V$
\begin{equation}\label{eq-bdry-color}
\begin{linkdiag}
\fill[gray!20] (0,0) rectangle (1.5,1);
\draw[thick,-o-={0.5}{<}] (0,0.4) -- node[pos=0.6,above]{$a$} (1.5,0.4);
\draw (1.5,0) -- +(0,1);
\fill (1.5,0.4) circle[radius=0.1] node[right]{$a$};
\end{linkdiag}
\quad\text{or equivalently}\quad
\begin{linkdiag}
\fill[gray!20] (0,0) rectangle (1.5,1);
\draw[thick,-o-={0.6}{>}] (0,0.4) -- node[pos=0.45,above]{$a^\ast$} 
(1.5,0.4);
\draw (1.5,0) -- +(0,1);
\fill (1.5,0.4) circle[radius=0.1] node[right]{$a$};
\end{linkdiag}.
\end{equation}
Note the conservation can take different forms using the above color
equivalence.  In diagrams, edge colors are often omitted for convenience.

\subsection{$\gsl$-skeins}
\label{sec-intro-skein}

We now introduce $\gsl$-skeins slightly informally, highlighting their main
features and use them to phrase our theorems, postponing the complete and
technical definitions to Section~\ref{sec-skein}.

The color lattice $\colorL$ determines the universal ring 
\begin{equation}
\label{RLambda}
\groundR=\ints[q^{\pm1},q^{\pm a_i},q^{\pm2a_ia_j}][\CC{b}^{-1},
b\in\colorL\setminus\{0\}], \qquad \CC{b}:=q^b-q^{-b}
\end{equation}
(where $\Lambda$ has $\BZ$-basis $1,a_1,\dotsc,a_r$) with 
field of fractions 
\begin{equation}
\groundK=\rats(q,q^{a_i},q^{2a_ia_j}) \,.
\end{equation}
Here, $q^\bullet$ are considered formal variables, but as usual in skein 
theory,
they can be specialized to values in a general ground ring by tensor product
(analogous to \cite[Prop.~2.2(4)]{Przytycki:fund} for the Kauffman bracket).

The $\colorL$-colored $\gsl$-skein module $\skein(M,V)$ (with $V$ omitted if
empty) is the quotient of the free $\groundR$-module generated by the isotopy
classes of framed and typically colored trivalent webs modulo the following
relations:
\begin{enumerate}
\item Split link:
\begin{equation}\label{eq-split-link}
\begin{linkdiag}
\fill[gray!20] (-0.2,0)rectangle(1,1);
\draw[thick] (0.5,0.5)circle(0.3) (0,0) -- (0,1);
\end{linkdiag}
=0.
\end{equation}
\item Bubble:
\begin{equation}\label{eq-bubble}
\begin{linkdiag}
\fill[gray!20] (-0.4,-0.3)rectangle(2,1.3);
\webnode{0.3,0.5}{L}{+};
\webnode{1.3,0.5}{R}{-};
\begin{scope}[thick]
\draw (-0.4,0.5)--(L);
\draw[-o-={0.5}{<}] (L)..controls (0.7,0.9) and (0.9,0.9).. node[pos=0.8,
above,inner sep=2pt]{$b$}(R);
\draw[-o-={0.5}{<}] (L)..controls (0.7,0.1) and (0.9,0.1).. node[pos=0.8,
below,inner sep=2pt]{$c$}(R);
\draw[-o-={0.8}{>}] (R) --node[pos=0.7,above]{$a$} (2,0.5);
\end{scope}
\end{linkdiag}
=
\frac{\CC{a+1}}{\CC{b}\CC{c}}
\begin{linkdiag}
\fill[gray!20] (0,0)rectangle(1.2,1);
\draw[thick,-o-={0.8}{>}] (0,0.5) --node[pos=0.8,above]{$a$} (1.2,0.5);
\end{linkdiag}.
\qquad
\begin{linkdiag}
\fill[gray!20] (-0.4,-0.3)rectangle(2,1.3);
\webnode{0.3,0.5}{L}{\nu};
\webnode{1.3,0.5}{R}{\nu};
\draw[thick] (-0.4,0.5)--(L) (R) -- (2,0.5)
  (L)..controls (0.7,0.9) and (0.9,0.9).. (R)
  (L)..controls (0.7,0.1) and (0.9,0.1).. (R);
\end{linkdiag}
=0.
\end{equation}
\item $I=H$:
\begin{equation}\label{eq-IH}
\begin{gathered}
\Iskein{+}{+}=\Hskein{+}{+}.\qquad
\Iskein{-}{-}=\Hskein{-}{-}.\\
\Iskein{+}{-}+\Iskein{-}{+}=\Hskein{+}{-}+\Hskein{-}{+}.
\end{gathered}
\end{equation}
\item Crossing:
\begin{equation}\label{eq-crossing}
\begin{linkdiag}
\fill[gray!20] (-0.3,-0.3)rectangle(1.3,1.3);
\begin{scope}[thick]
\begin{knot}
\strand (1.2,1.3) -- (-0.2,-0.3);
\strand (-0.2,1.3) -- (1.2,-0.3);
\end{knot}
\path[-o-={0.8}{>}] (1.2,1.3) -- node[pos=0.7,left]{$a\vphantom{b}$} 
(-0.2,-0.3);
\path[-o-={0.8}{>}] (-0.2,1.3) -- node[pos=0.7,right]{$b$} (1.2,-0.3);
\end{scope}
\end{linkdiag}
=q^{-2ab-a-b}\left(\CC{a-b}\Hskein{+}{-}-q^b\CC{a}\Iskein{-}{+}
+q^{-a-1}\CC{b+1}\Iskein{+}{-}\right).
\end{equation}
\item Perturbation: Suppose $\alpha_i$ are webs identical outside a ball
$B^3\subset M$, and $D_i=\alpha_i\cap B^3$ have perturbations 
$D_i(\varepsilon)$
satisfying a relation
\begin{equation}
\sum_i c_i(\varepsilon) D_i(\varepsilon)=0.
\end{equation}
If all $c_i(\varepsilon)$ are finite at $\varepsilon=0$, then we also 
impose the
relation for the webs at $\varepsilon=0$.
\end{enumerate}
This last relation is defined formally using a perturbation
$\colorL_\varepsilon$ of the lattice and the corresponding color ring; see
Equation~\eqref{lambdae} below.

Note each skein relation also carries the implicit condition that all edges have
typical colors. The most notable one being the crossing relation
\eqref{eq-crossing} where the specialization $a=b$ would make some colors
atypical and the relation as stated does not apply.

There are several consequences of the above skein relations (including 
how to resolve a negative crossing, and how to resolve crossings labeled by 
the same color) discussed in Lemma~\ref{lemma-more-rel} below.

The next theorem provides a diagrammatic representation of invariant 
tensors via
a map defined by the usual Reshetikhin--Turaev construction.

Fix a list $(a_0,a_1,\dotsc, a_n)$ for $n\ge1$ and consider 
$V_{n+1}\subset\partial B^3$ is colored by $(a_0,a_1,\dotsc,a_n)$. 

\begin{theorem}
\label{thm.1}
\begin{enumcc}[(1)]
\item The skein module $\skein(B^3,V_{n+1})$ is a free $\groundR$-module of rank 
\begin{equation}
\label{m0nl}
m^{(n)}_{0,\ell} = 
\frac{1}{n}\binom{n}{\ell}\binom{n}{\ell+1}
\end{equation}
if $\ell:=-1-\displaystyle\sum_{i=0}^na_i \in \{0,1,\dotsc,n-1\}$, and zero
otherwise.
\item Over the field $\groundK$, 
\begin{equation}
\label{inv0nl}
\dim_{\groundK} 
\Inv_{\Uqsl}\Big(\bigotimes_{i=0}^n V(0,a_i)\Big) = m^{(n)}_{0,\ell} \,.
\end{equation}
\item There are explicit bases for the domain and codomain
and a $\groundK$-vector space isomorphism
\begin{equation}
\label{eq-Fbar}
F:\skein(B^3,V_{n+1})\otimes_{\groundR}\groundK\to\Inv_{\Uqsl}\Big
(\bigotimes_{i=0}^n V(0,a_i)\Big)
\end{equation}
induced by the Reshetikhin--Turaev construction.
\end{enumcc}
\end{theorem}

Part (2) follows from Proposition~\ref{prop.mult}. Parts (1) and (3) are 
proved in Section~\ref{sec-proof}. Note that the theorem holds without
any genericity assumption on the colors $a_i$.

\begin{remark}
We briefly discuss the two remaining skein modules for $B^3$.
\begin{enumcc}[(1)]
\item If there is no boundary point, then
$\skein(B^3,\emptyset)=\skein(S^3)=\groundR$ as a module. The value of a link
with proper normalization is the Links--Gould polynomial. See
Proposition~\ref{prop-S3}.
\item If there is exactly one boundary point, then it is easy to show that both
the skein module $\skein(B^3,V_1)$ and the invariant space
$\Inv_{\Uqsl}(V(0,a))$ are zero for typical colors.
\end{enumcc}
\end{remark}

\begin{example}\label{exa-basis1}
When the list of colors excluding $a_0$ is generic (as defined in
Section~\ref{sub.mult}), a basis of the skein module $\skein(B^3,V_6)$ for $n=5$
and $\ell=2$ consists of $m^{(5)}_{2,0}=20$ trees shown in
Figure~\ref{fig-basis-trees}. Since these webs are trees, the edge colors are
uniquely determined by the labeling $\pm$ of their vertices.

\begin{figure}
\centering
\includegraphics[width=0.8\linewidth]{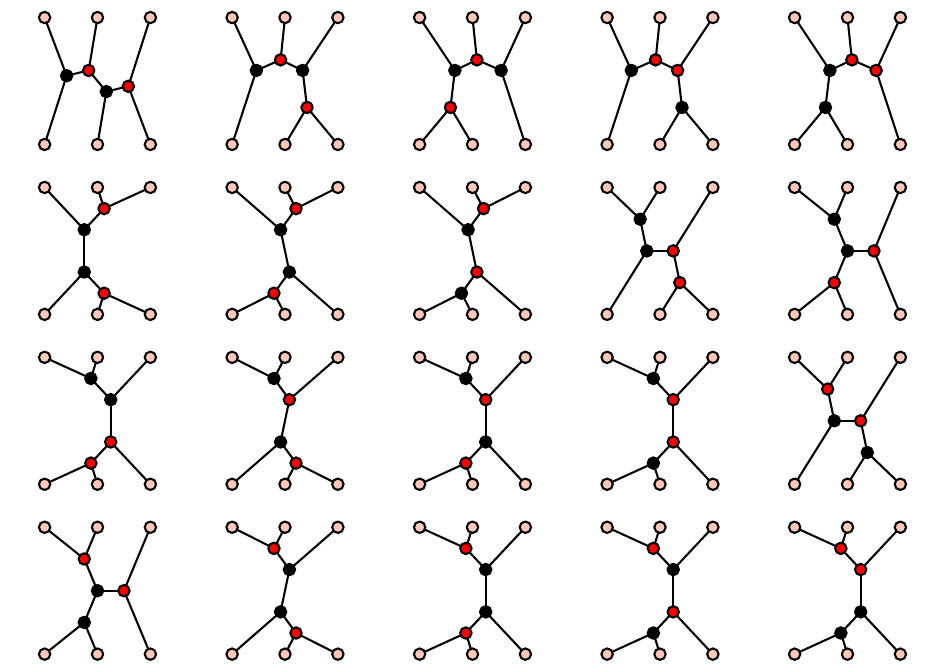}
\caption{A basis of 20 trees in $\skein(B^3,V_6)$ with generic colors and
$\ell=2$, with black and red vertices are used in place of $+,-$ vertices,
respectively.}
\label{fig-basis-trees}
\end{figure}
\end{example}

In turn, this implies the following diagrammatic description of the admissible
skein module $\skein_\mathcal{I}(M,V)$ by \cite{CGP:admissible}, where
$\mathcal{I}$ denotes the ideal of projective objects in $\Repsl$ (a precise
definition is given in Section~\ref{sec.sl21} below).

\begin{theorem}
\label{thm.2}
For every oriented 3-manifold $M$ and $V \subset \partial M$ a finite set of
framed points colored by $\colorL\setminus\{0,-1\}$, we have an isomorphism 
of
vector spaces
\begin{equation}
\skein(M,V)\otimes_{\groundR}\groundK \simeq \skein_\mathcal{I}(M,V) \,.
\end{equation}
\end{theorem}

In the next remark we point out some features of $\gsl$-webs that are 
absent from their relative webs of simple Lie algebras.

\begin{remark}
\label{rem-features}
\begin{enumcc}[(1)]
\item When $M=\Sigma\times[0,1]$ and $V=\emptyset$, the $\gsl$-skein module has
a product by stacking in the $[0,1]$-direction as usual, but the product lacks a
unit since the empty set is not an element of the skein module. If it were, the
split link relation \eqref{eq-split-link} would kill a lot of elements.
\item The specialization $q=1$, where one expects to obtain the coordinate ring
of the classical character variety, is not possible since the ring $\groundR$
contains $\frac{1}{\{1\}}=\frac{1}{q-q^{-1}}$. One can specialize $q=e^\hbar$
resulting in the field $\BQ(q^{a_i},q^{2a_ia_j})(\!(\hbar)\!)$. Nevertheless,
when $M=S^3$, the normalized evaluation of every colored link gives an element
of $\ints[q^{\pm1},q^{\pm a_i}]$ (the colored Links--Gould polynomial), which
can be specialized at $q=1$, and then equals to the square of the multivariable
Alexander polynomial; see Theorem~\ref{thm-alex} below.
\end{enumcc}
\end{remark}

\begin{question}
\label{que.q=1}
What is the classical limit $q=1$ of $\skein(M,V)$?
\end{question}

\begin{remark}\label{rem-stacking}
When $M=[0,1]^3$ and $V\subset\{\frac{1}{2}\}\times[0,1]\times\{0,1\}$, the
skein modules $\skein([0,1]^3,V)$ has two stackings, which are not to be
confused with Remark~\ref{rem-features}(1), making them into a braided monoidal
(in fact ribbon) category. The composition of morphisms is given by stacking in
the third direction. The tensor product of morphisms is given by stacking in the
second direction. Braiding is given by~\eqref{eq-crossing}.
\end{remark}

\subsection{Skeins and $\gsl$-Hecke algebras}
\label{sub.Bn}

In this section we explain how our $\gsl$ webs are related to some key results
of Marin--Wagner on the structure of cubic Hecke algebras and their
quotients~\cite{Marin-Wagner}, and to a theorem of Anghel~\cite{Anghel}. To do
so, we choose $\colorL=\ints+a\ints$ and specialize to the skein module
$\skein(B^3,V_{n,n})$ of the 3-disk $B^3$ and $V_{n,n}$ consisting of $n$ points
colored $a$ and $n$ points colored $a^\ast$. Using a writhe correction similar
to the Jones polynomial from the Kauffman bracket, we can eliminate the
generator $q^{2a^2}$ from coefficients. (See e.g.\ \eqref{eq-ML}.) Therefore, we
work with the coefficient ring $R_a$ and its field of fractions $k_a$
\begin{equation}
R_a=\ints[q^{\pm1},q^{\pm a}][\CC{b}^{-1},b\in\colorL\setminus\{0\}],
\quad
k_a=\rats(q,q^a).
\end{equation}

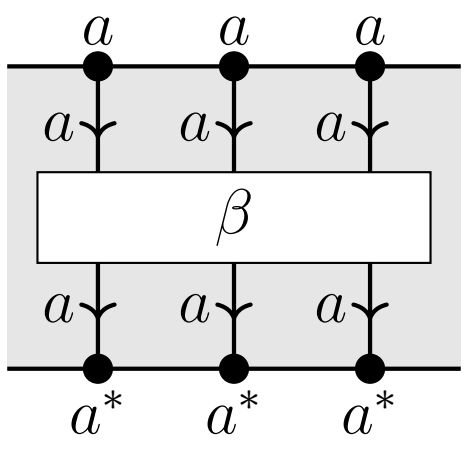
\begin{figure}
\centering
\begin{tikzpicture}
\fill[gray!20] (0,0) rectangle (3,2);
\draw[thick] (0,0) -- (3,0) (0,2) -- (3,2);
\foreach \x in {0.6, 1.5, 2.4} {
\fill (\x,2) circle[radius=0.1] node[above]{$a$};
\draw[thick,-o-={0.6}{>}] (\x,2) -- node[left]{$a$} (\x,1.2);
\fill (\x,0) circle[radius=0.1] node[below]{$a^\ast$};
\draw[thick,-o-={0.6}{>}] (\x,0.8) -- node[left]{$a$} (\x,0);
}
\draw[fill=white] (0.2,0.7) rectangle (2.8,1.3) (1.5,1)node{$\beta$};
\end{tikzpicture}
\caption{Elements of $\skein(B^3,V_{n,n})$.}\label{fig-hecke}
\end{figure}

By convention \eqref{eq-bdry-color}, elements in the skein module
$\skein(B^3,V_{n,n})$ look like Figure~\ref{fig-hecke}. This skein module is an
associative $R_a$-algebra by composition, as explained in
Remark~\ref{rem-stacking}, and it is a cousin of the famous
Temperley--Lieb algebra. Note the algebra has a unit given by the identity braid
in Figure~\ref{fig-hecke}.

The next corollary of Theorem~\ref{thm.1} gives a diagrammatic 
presentation of the invariant tensors of the tensor powers of $V(0,a)$.

\begin{corollary}
$\skein(B^3,V_{n,n})$ is a free $R_a$-module of rank
\begin{equation}
  \label{m0nn}
m^{(2n-1)}_{0,n-1}
=\frac{1}{2n-1}\binom{2n-1}{n-1}\binom{2n-1}{n}
=\frac{(2n-2)!(2n-1)!}{\big((n-1)!n!\big)^2}
\end{equation}
and the Reshetikhin--Turaev map
\begin{equation}
\label{2bsn}
h_n: \skein(B^3,V_{n,n})\otimes_{R_a}k_a 
\to \End_{\Uqsl}(V(0,a)^{\otimes n})
\end{equation}
is a $k_a$-algebra isomorphism.
\end{corollary}
The dimension count is a new proof of a theorem of 
Anghel~\cite[Thm.~1]{Anghel}.

\begin{example}\label{exa-basis2}
Note the coloring list $(a,\dotsc,a,a^\ast,\dotsc,a^\ast)$ of the boundary
points is not generic, so the basis is different from Example~\ref{exa-basis1}.
Here, a basis of the skein module $\skein(B^3,V_{n,n})$ consists of 20 graphs
shown in Figure~\ref{fig-basis-end}. These graphs are obtained from
Figure~\ref{fig-basis-trees} by attaching some fixed tree which perturbs the
colors to be generic. See Section~\ref{sec-proof}.

\begin{figure}
\centering
\includegraphics[width=0.8\linewidth]{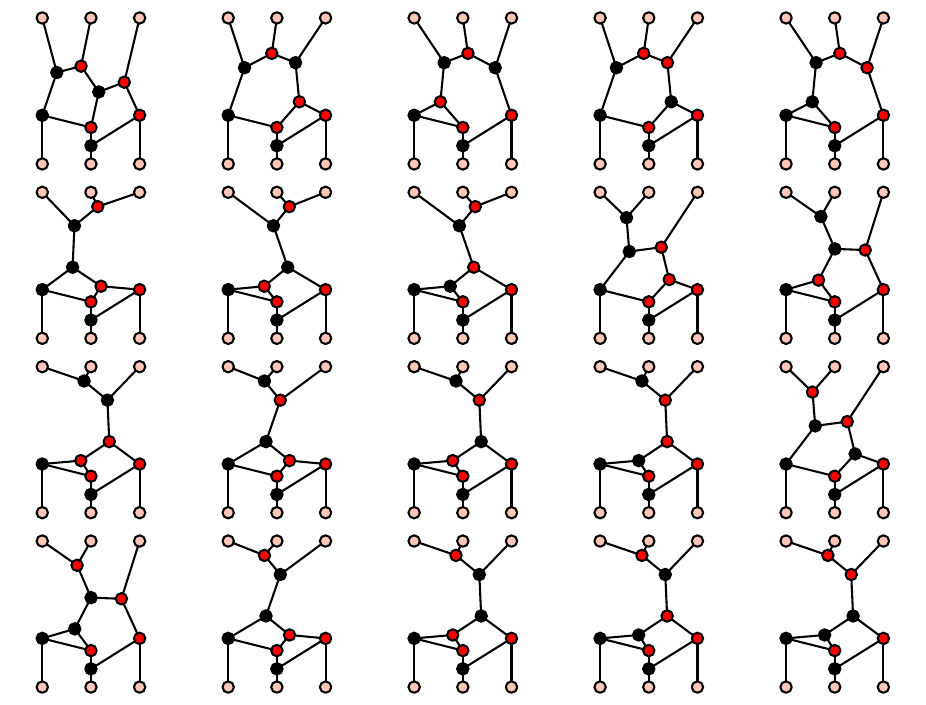}
\caption{20 basis graphs in $\skein(B^3,V_{n,n})$.}
\label{fig-basis-end}
\end{figure}
\end{example}

The formula resolving a crossing and the very definition of the skein 
module whose elements consist of isotopy classes of webs
gives a compatible collection of representations of the braid
group $B_n$ on $n$ strands
\begin{equation}
\label{rhoBn}
\rho_n: R_a[B_n] \to \skein(B^3,V_{n,n}) \,.
\end{equation}
This can be used to compute the colored
Links--Gould polynomials of links in $S^3$.
In a subsequent publication we will show that the map~\eqref{rhoBn}
is onto, which combined with map~\eqref{2bsn} implies that
the composition
\begin{equation}
k_a[B_n] \twoheadrightarrow \skein(B^3,V_{n,n}) \otimes_{R_a} k_a 
\cong \End_{\Uqsl}(V(0,a)^{\otimes n})
\end{equation}
is onto, recovering a theorem of Marin--Wagner~\cite[Thm.~0.1]{Marin-Wagner}.

A second theorem of Marin--Wagner~\cite{Marin-Wagner} concerns the 
quotient 
\begin{equation}
\label{Cndef}
A_n=k_a[B_n]/(R_1,R_2,R_3),
\end{equation}
of the group algebra $k_a[B_n]$ by the 2-sided ideal $(R_1,R_2,R_3)$ given
explicitly in \cite{Ga:skein}, where the variables
$(t_0,t_1)=(q^{2a},q^{-2a-2})$. $A_n$ is a quotient of the cubic Hecke algebra
$k_a[B_n]/(R_1)$. In \cite{Marin-Wagner} (see also~\cite{Ga:skein}) it was shown
that the map $k_a[B_n] \to \End_{\Uqsl}(V(0,a)^{\otimes n})$ descends to a map 
\begin{equation}
\label{LGn}
\rho_n^{\LG}: A_n \to \End_{\Uqsl}(V(0,a)^{\otimes n}).
\end{equation}

One can check that the 3, 12 and 80-term relations
$R_1$, $R_2$ and $R_3$ on the braid groups $B_2$, $B_3$ and $B_4$ 
respectively, are mapped to 0 under $\rho_n \otimes k_a$. 
We will discuss this in detail in a subsequent publication.

Putting everything together, we have a commutative triangle
\begin{equation}
\label{SkC}  
\begin{tikzcd}[row sep=tiny]
& \skein(B^3,V_{n,n})\otimes_{R_a}k_a \arrow[dd,"h_n","\cong"']\\
A_n \arrow[ru,two heads,"\rho_n \otimes k_a"] 
\arrow[rd,two heads,"\rho_n^{\LG}"']\\
& \End_{\Uqsl}(V(0,a)^{\otimes n})
\end{tikzcd}
\qquad
\begin{tikzcd}[row sep=tiny]
& \text{$\gsl$-Webs} \arrow[dd,"h_n","\cong"']\\
\text{Braids} \arrow[ru,two heads,"\rho_n \otimes k_a"] 
\arrow[rd,two heads,"\rho_n^{\LG}"']\\
& \text{$\gsl$-Tensors}
\end{tikzcd}
\end{equation}
The map $\rho_n \otimes k_a$ from braids to $\gsl$-webs 
gives a pictorial realization of the map $\rho_n^{\LG}$ from braids to
$\gsl$-invariant tensors, answering a question of Marin--Wagner.

\subsection{Cabling and the $\gsl$-skein algebra of the annulus}
\label{sub.cabling}

In this section we study the $\gsl$-skein module 
$\skein(A):=\skein(A\times[0,1],\emptyset)$ of the annulus $A$, which
is a commutative algebra with the stacking product.
By the standard argument, skein modules have the following
functorial property.

\begin{lemma}
\label{lem-ann-embed} 
Suppose $\alpha$ is an $\gsl$-web in $(M,V)$, and $K\subset M$ is a framed knot
disjoint with $\alpha$. Let $i:A\times[0,1]\embed M$ denote the inclusion of the
solid torus as a regular neighborhood of $K$ with $A$ aligned with the framing
of $K$. Then $i$ induces a homomorphism of $\groundR$-modules
\begin{equation}
i_\ast:\skein(A)\to\skein(M,V),\qquad
i_\ast(\beta)=\alpha\sqcup i(\beta)
\end{equation}
for every $\gsl$-web $\beta$ in $A\times[0,1]$.
\end{lemma}

As we will explain, the Links--Gould polynomials are just normalized skeins in $S^3$.
Therefore, to obtain cabling (or more generally, satellite) formulas, we can
calculate the skein of the pattern in the skein algebra $\skein(A)$ and then
combine with the rest by functoriality. We apply this principle to $(p,0)$- and
$(p,1)$-cables, given below in Theorems~\ref{thm-ML-cable} and
\ref{thm-ML-recursion}.

To phrase our results, consider an oriented link $L$ in $S^3$ with
$r$ ordered components and a list of positive integers $n_1,\dotsc,n_r$. Then, one can define a polynomial invariant 
\begin{equation}
\label{LGmc1}
\LG^{(n_1,\dotsc,n_r)}_L(q;t_1,\dotsc,t_r) 
\in \BZ[q^{\pm 1},t_1^{\pm 1},\dotsc,t_r^{\pm 1}], \qquad
\LG^{(n)}_{\text{Unknot}}(q;t) = 1 
\end{equation}
following \cite{GP:multi}. We discuss its definition 
in Section~\ref{sub.LGdef} below, along with its useful renormalization
$\GG$ that satisfies the symmetries 
\begin{equation}
\label{LGmc2}
\GG^{(n_1,\dotsc,n_r)}_{\sigma_i L}(q;t_1,\dotsc,t_r) =
\GG^{(n_1,\dotsc,n_r)}_{L}(q;t_1,\dotsc,t_i^\ast,\dotsc,t_r),
\end{equation}
where $t_i^\ast = q^{-n_i} t_i^{-1}$ and $\sigma_i L$ is obtained from $L$ by
reversing the orientation of its $i$-th component, and
\begin{equation}
\label{LGmc3}
\begin{aligned}
\GG^{(n_1,\dotsc,n_r)}_{L}(q;t_1^\ast,\dotsc,t_r^\ast) &=
\GG^{(n_1,\dotsc,n_r)}_{L}(q;t_1,\dotsc,t_r) 
\\
\GG^{(n_1,\dotsc,n_r)}_{\text{mirror}(L)}(q^{-1};t_1,\dotsc,t_r) &=
\GG^{(n_1,\dotsc,n_r)}_{L}(q^{-1};t_1,\dotsc,t_r) 
\,.
\end{aligned}
\end{equation}

When $n_i=1$ and $t_i=t$ for all $i$  this is the original Links--Gould
polynomial~\cite{Links-Gould}. When $n_i=1$ for all $i$, we suppress them from
the notation.

Fix an oriented link $L\subset S^3$ with $r$ ordered components, and let
$\mu_{i,j}$ denote its linking matrix. The next theorem proves a conjecture of
\cite[Remark~4.1]{Ga:colored} and extends it to links. The coefficients in this
result and the next are slightly different from \cite{Ga:colored} due to normalizations.

\begin{theorem}
\label{thm-ML-cable}
With $L$ as above, let $L(n,1)$ be the link $L$ with the last component
replaced by its $(n,1)$-cable. Then
\begin{multline}
(-1)^{n-1}\GG^{(n_1,\dotsc,n_{r-1},1)}_{L(n,1)}(q;t_1,\dotsc,t_r)
=\GG^{(n)}_L(q;t_1,\dotsc,t_{r-1},t_r^n)\\
+\sum_{i=1}^{n-1}(-1)^i\big(
Q_1^i\GG^{(n-i)}_L(q;t_1,\dotsc,t_{r-1},t_r^n)
+q^{-2i}Q_1^{-i}\GG^{(n-i)}_L(q;t_1,\dotsc,t_{r-1},q^it_r^n)\big),
\end{multline}
where $\GG^{(n)}_L$ is the shorthand of $\GG^{(n_1,\dotsc,n_{r-1},n)}_L$
and $Q_1=t_r^2\prod_{i=1}^{r-1}(q^{n_i}t_i^2)^{\mu_{r,i}}$.
\end{theorem}

The next theorem extends a cabling formula of \cite[Thm.1.1]{Ga:colored}
to links. Again, fix an oriented link $L\subset S^3$ with $r$ ordered 
components.

\begin{theorem}
\label{thm-ML-recursion}
With $L$ as above, let $L(p+1,0)$ be the link obtained by adding $p$
parallel copies of the last component as the $(r+1)$-th to the $(r+p)$-th
components. Then the normalized invariants satisfy a linear relation
\begin{multline}
\label{eq-ML-recursion}
\GG^{(n_1,\dotsc,n_{r-1},1,\dotsc,1)}_{L(p+1,0)}(q;t_1,\dotsc,t_{r-1},t_r,
\dotsc,t_{r+p})\\
=\sum_{k=0}^p\sum_{\ell=0}^{p-k}m^{(p+1)}_{k,\ell}Q^{p-k-2\ell}
\GG^{(n_1,\dotsc,n_{r-1},k+1)}_L(q;t_1,\dotsc,t_{r-1},q^\ell\prod_{i=r}^{r+p}
t_i), \quad Q=\prod_{i=1}^{r-1}(q^{n_i}t_i^2)^{\mu_{r,i}}
\end{multline}
where $m^{(n)}_{k,\ell}$ is the multiplicity defined in \eqref{otimesn}
and computed in~\eqref{eq-mult}.
\end{theorem}

Note that there is exactly one term with the highest labels $(n_1,\dotsc,n_r)$
in the lexicographical order at $(k,\ell)=(p,0)$, and this can be used
recursively to calculate any invariant $\GG^{(n_1,\dotsc,n_r)}_L$ from the basic
ones $\GG_L$.

\subsection{The colored link invariants at $q=1$}
\label{sub.GLq=1}

Let $\Delta_L(t_1,\dotsc,t_r) \in \BZ[t_1^{\pm 1},\dotsc,t_r^{\pm 1}]$ 
denote the multivariable Alexander polynomial~\cite{Tu:alex} of an oriented link $L$ in $S^3$ with $r$ ordered components. 
A careful definition of this polynomial (which is often defined only
up to multiplication by $\pm 1$ and by a monomial in the variables $t_i$)
is given in~\cite{Boyer}, \cite{Jiang} (where a braid-type skein theory
with a unique solution and a reduction algorithm is also given), as well
as \cite{Cimasoni} and \cite{Conway:pot}. 

A useful normalization of the multivariable Alexander polynomial
is the Conway potential function defined by
\begin{equation}
\label{na}  
\nabla_L(t_1,\dotsc,t_r)=\begin{cases}
\Delta_L(t_1^2)/(t_1-t_1^{-1}),&r=1,\\
\Delta_L(t_1^2,\dotsc,t_r^2),&r\ge2.
\end{cases}
\end{equation}
Its diagonal specialization satisfies the Conway skein relation
\begin{equation}
\nabla_{L^+}(t,\dotsc,t) - \nabla_{L^-}(t,\dotsc,t) = (t-t^{-1})
\nabla_{L^0}(t,\dotsc,t),
\end{equation}
and the Conway potential function satisfies the symmetries
\begin{equation}
\nabla_{\sigma_i L}(t_1,\dotsc,t_r) = - \nabla_{L}(t_1,\dotsc,
t_i^{-1}, \dotsc, t_r)
\end{equation}
where $\sigma_i L$ is the link obtained from $L$ by reversing the 
orientation of its $i$-th component, and
\begin{equation}
\begin{split}
\nabla_L(t_1^{-1},\dotsc,t_r^{-1}) &= (-1)^r \nabla_L(t_1,\dotsc,t_r),\\
\nabla_{\text{mirror}(L)}(t_1,\dotsc,t_r)
&= (-1)^{r-1}\nabla_L(t_1,\dotsc,t_r)
\,.
\end{split}
\end{equation}

\begin{theorem}\label{thm-alex}
For every oriented link $L\subset S^3$ with $r$ ordered components, 
we have
\begin{equation}
\label{eq-higher-alex}
\GG^{(n_1,\dotsc,n_r)}_L(1;t_1,\dotsc,t_r)=
\Big(\prod_{i=1}^r(-1)^{n_i-1}n_i\Big)\nabla_L(t_1,\dotsc,t_r)^2.
\end{equation}
\end{theorem}

Note that~\eqref{eq-higher-alex} is equivalent to
\begin{equation}
\label{LGalex}
\LG^{(n_1,\dotsc,n_r)}_L(1;t_1,\dotsc,t_r)=
\begin{cases}
\Delta_L(t_1^2)^2,&r=1,\\
\Big(\prod_{i=1}^r(-1)^{n_i-1}n_i\Big)\Delta_L(t_1^2,\dotsc,t_r^2)^2,&r\ge2.
\end{cases}
\end{equation}

\subsection*{Acknowledgements}

The authors wish to thank Renaud Detcherry and Roland van der Veen
for enlightening conversations.

\noindent
{\bf AI use disclosure.}
No LLMs were used in the preparation and writing of this paper.

%%%%%%%%%%%%%%%%%%%%%%%%%%%%%%%%%%%%%%%%%%%%%%%%%%%%%%%%%%%%%%%%%%%%%%%%%%%%
%%%%%%%%%%%%%%%%%%%%%%%%%%%%%%%%%%%%%%%%%%%%%%%%%%%%%%%%%%%%%%%%%%%%%%%%%%%%

\section{$\Uqsl$ and its representation theory}
\label{sec.sl21}

In this section we recall the Hopf superalgebra $\Uqsl$ and its representation
category $\Repsl$. All results in this section are over the field $\groundK$.
Recall the supercommutator $[x,y]=xy-(-1)^{d(x)d(y)}yx$.

\begin{definition}
$\Uqsl$ is the Hopf superalgebra generated by even elements $K_1^{\pm1}$,
$K_2^{\pm1}$, $E_1,F_1$ and odd elements $E_2,F_2$, satisfying the following
relations.
\begin{equation}
\begin{aligned}
K_iE_j&=q^{a_{ij}}E_jK_i,&
E_1E_{12}&=qE_{12}E_1,&
&E_2^2=F_2^2=0,\quad
[K_i,K_j]=0,
\\
K_iF_j&=q^{-a_{ij}}F_jK_i,&
F_1F_{12}&=qF_{12}F_1.&
&[E_i,F_j]=\delta_{ij}\frac{K_i-K_i^{-1}}{q-q^{-1}},
\end{aligned}
\end{equation}
where
\begin{equation}\label{eq-Cartan}
(a_{ij})=\begin{pmatrix}2&-1\\-1&0\end{pmatrix}
\end{equation}
is the Cartan matrix, and
\begin{equation}
E_{12}=E_1E_2-q^{-1}E_2E_1,\qquad
F_{12}=F_2F_1-qF_1F_2.
\end{equation}
The coproduct, counit, and antipode are given by
\begin{equation}
\begin{aligned}
\Delta(E_i)&=E_i\otimes1+K_i^{-1}\otimes E_i,&
\epsilon(E_i)&=0,&
S(E_i)&=-K_iE_i,\\
\Delta(F_i)&=F_i\otimes K_i+1\otimes F_i,&
\epsilon(F_i)&=0,&
S(F_i)&=-F_iK_i^{-1},\\
\Delta(K_i)&=K_i\otimes K_i,&
\epsilon(K_i)&=1,&
S(K_i)&=K_i^{-1}.\\
\end{aligned}
\end{equation}
Let $H_1,H_2\in\gsl$ be the Cartan elements formally related to
$K_1,K_2$ as $K_i=q^{H_i}$.
\end{definition}

\subsection{Roots and positivity}
\label{sec-roots}

Let $\mathfrak{h}=\Span\{H_1,H_2\}$. Weights and roots are elements of
$\mathfrak{h}^\ast$ by definition, and we parametrize them using their 
values on $H_1,H_2$.

\begin{figure}
\centering
\begin{tikzpicture}
\draw[gray] (-2.5,-1.5) grid (2.5,1.8);
\draw[->] (-2.5,0) -- (2.5,0) node[right]{$H_1$};
\draw[->] (0,-1.5) -- (0,1.8) node[above]{$H_2$};
\begin{scope}[thick, inner sep=1pt]
\draw[->] (0,0) -- (2,-1)node[below left]{$\alpha_1$};
\draw[->] (0,0) -- (1,-1)node[below left]{$\alpha_{12}$};
\draw[->] (0,0) -- (-1,0)node[below left]{$\alpha_2$};
\draw[->] (0,0) -- (-2,1)node[above right]{$-\alpha_1$};
\draw[->] (0,0) -- (-1,1)node[above right]{$-\alpha_{12}$};
\draw[->] (0,0) -- (1,0)node[above right]{$-\alpha_2$};
\end{scope}
\fill[blue,opacity=0.2] (-2.5,0) -- (0,0) -- (2.5,-1.25) |- (-2.5,-1.5) -- cycle;
\fill[red,pattern=north east lines,opacity=0.5] (0,0) rectangle (2.5,-1.5);
\end{tikzpicture}
\caption{Roots and positivity.}
\label{fig-roots}
\end{figure}
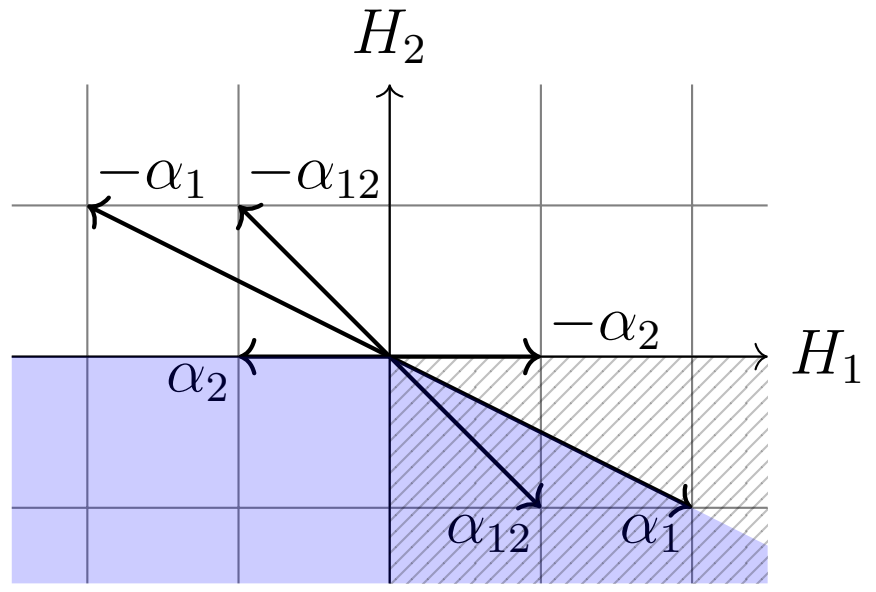

$\gsl$ has six roots $\{\alpha_1,\alpha_2,\alpha_{12},
-\alpha_1,-\alpha_2,-\alpha_{12}\}$ shown in Figure~\ref{fig-roots}
where $\alpha_{12}=\alpha_1+\alpha_2$, with
$\alpha_1$ even and $\alpha_2$ and $\alpha_{12}$ odd.
The standard choice of positivity, used for the
Cartan matrix \eqref{eq-Cartan}, declares the positive roots to be
$\Phi^+=\{\alpha_1,\alpha_2,\alpha_{12}\}$.  
The generators $E_\ast$ correspond to the roots
$\alpha_\ast$, and $F_\ast$ correspond to $-\alpha_\ast$. 
The positive cone is the blue region in Figure~\ref{fig-roots}.

There is an alternative notion of positivity where $-\alpha_2,\alpha_{12}$ are
the simple roots. This choice appears as the 5th entry in \cite[Table~V]{Kac}.
We extend the cone of positivity to be spanned by
$-\alpha_2,\alpha_1+2\alpha_2$, which is the shaded lower right quadrant in
Figure~\ref{fig-roots}. The reason of this alternative choice of positivity and
the corresponding cone is for Lemma~\ref{lemma-tree-high}. See also
Remark~\ref{rem-basis}.

Intuitively, given a set of weights plotted on the $(H_1,H_2)$-plane (such as
Figure~\ref{fig-basis-weights}), the standard highest weight is the bottom left
corner, while the alternative highest weight is the bottom right.

In particular, for the typical representations $V(k,a)$ discussed below, if $v$
is the highest weight vector in the standard choice of positivity, then $F_2v$
has the highest weight in the alternative order. Conversely, if $v_1$ has the
highest weight in the alternative order, then $E_2v_1$ has the highest weight in
the standard order.

\subsection{The category of representations}

A $\Uqsl$-weight module is a finite dimensional $\Uqsl$-module such that the
Cartan generators $K_i$ acts diagonally on a basis of weight vectors. Let
$\Repsl$ denote the category of $\Uqsl$-weight modules. It has the structure of
a ribbon category. See e.g. \cite{Ga:skein}. Here, we describe the family
$V(k,a)$ of representations.

$V(k,a)$ is a $4(k+1)$-dimensional space deformed from Kac modules \cite{Kac}.
It is generated by a vector $v$ such that
\begin{equation}
K_1v=q^kv,\qquad
K_2v=q^av,\qquad
F_1^kv\ne0,\qquad
F_1^{k+1}v=0.
\end{equation}
$v$ is the highest weight vector in the usual order, but $F_2v$ has the highest
weight in the alternative order. The vector $v$ can be even or odd. In this
paper, $v$ has parity $k\bmod{2}$ for $V(k,a)$.

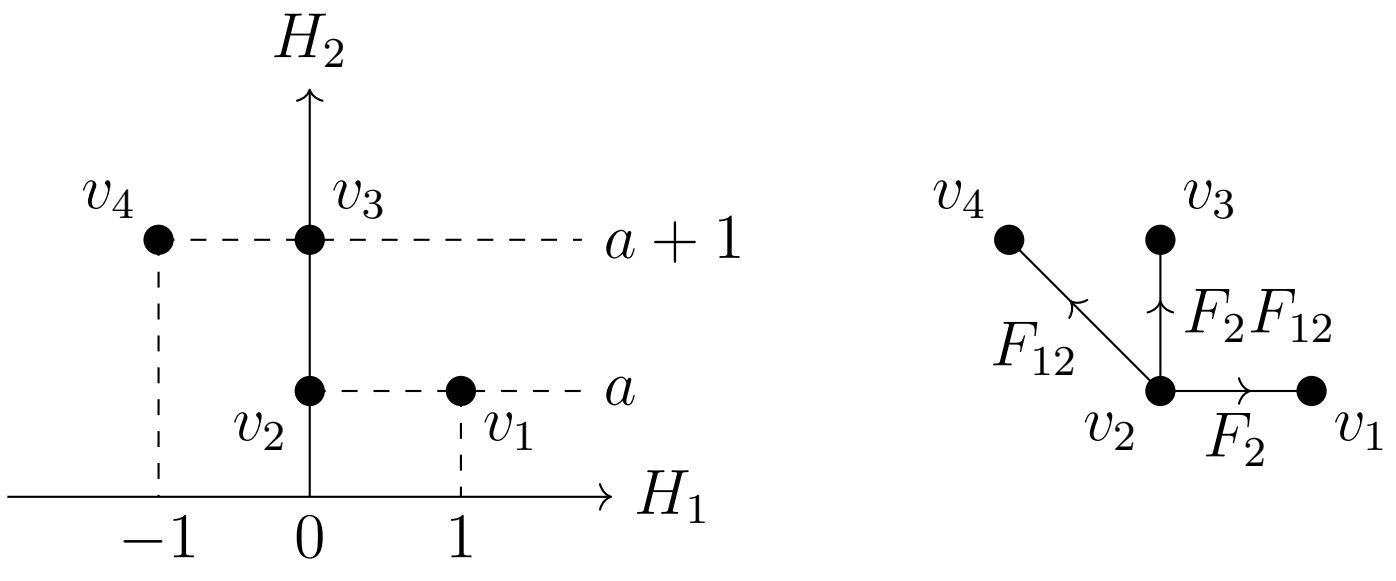
\begin{figure}
\centering
\begin{tikzpicture}[baseline=0cm]
\draw[->] (-2,-0.2) -- (2,-0.2) node[right]{$H_1$};
\draw[->] (0,-0.2) -- (0,2.5) node[above]{$H_2$};
\fill[radius=0.1] (0,0.5)circle[] node[below left]{$v_2$}
  (1,0.5)circle[] node[below right]{$v_1$}
  (0,1.5)circle[] node[above right]{$v_3$}
  (-1,1.5)circle[] node[above left]{$v_4$};
\draw[dashed,anchor=base] (1,0.5) -- (1,-0.2) +(0,-0.4)node{$1$}
  (0,-0.2) +(0,-0.4)node{$0$} (-1,1.5) -- (-1,-0.2) +(0,-0.4)node{$-1$};
\draw[dashed] (-1,1.5) -- (1.8,1.5) node[right]{$a+1$}
  (0,0.5) -- (1.8,0.5) node[right]{$a$};
\end{tikzpicture}
\qquad
\begin{tikzpicture}[baseline=0cm]
\fill[radius=0.1] (0,0.5)circle[] node[below left]{$v_2$}
  (1,0.5)circle[] node[below right]{$v_1$}
  (0,1.5)circle[] node[above right]{$v_3$}
  (-1,1.5)circle[] node[above left]{$v_4$};
\draw[-o-={0.6}{>}] (0,0.5) -- node[below]{$F_2$} (1,0.5);
\draw[-o-={0.6}{>}] (0,0.5) -- node[right]{$F_2F_{12}$} (0,1.5);
\draw[-o-={0.6}{>}] (0,0.5) -- node[below left,inner sep=1pt]{$F_{12}$} (-1,
1.5);
\end{tikzpicture}
\caption{Weights of the standard basis of $V(0,a)$.}
\label{fig-basis-weights}
\end{figure}

$V(k,a)$ has a standard PBW basis
\begin{equation}
v_{i,1}=F_2F_1^iv,\quad v_{i,2}=F_1^iv,\quad
v_{i,3}=F_2F_{12}F_1^iv,\quad v_{i,4}=F_{12}F_1^iv,\qquad
i=0,\dotsc,k.
\end{equation}
For $k=0$, the first subscript $i=0$ is omitted for brevity. The weights of the
basis vectors of $V(0,a)$ is shown in Figure~\ref{fig-basis-weights}.

\begin{remark}\label{rem-basis}
The basis of $V(0,a)$ in \cite{GP:multi} is
\begin{equation}
v'_1=q^av_2,\quad v'_2=-q^{-1}v_4,\quad
v'_3=-q^{-1}v_3,\quad v'_4=v_1.
\end{equation}
Yet another basis, used in \cite{Ga:colored}, is
\begin{equation}
v''_1=v=v_2,\quad v''_2=F_2v''_1=v_1,\quad
v''_3=F_1v''_2=-q^{-1}v_4,\quad v''_4=F_2v''_3=-q^{-1}v_3.
\end{equation}
Since all these bases consist of weight vectors, they only differ in order and
scaling. The basis $v''_i$ is in the usual weight order, as demonstrated by the
actions of $F_1,F_2$. We choose a different order for
Lemma~\ref{lemma-tree-high}. In fact, we only require $v_1$ to be the highest
and $v_4$ to be the lowest, but we made a choice of order between $v_2$ and
$v_3$ so that the $+$-vertex ranks higher than the $-$-vertex in
Section~\ref{sec-sk-highest}. The scalings are inconsequential in the proofs but
convenient for calculations by programs.
\end{remark}

$V(k,a)$ is called typical if $a\ne0,-k-1$.

\begin{lemma}[{\cite[Lemmas~3.4 and 3.5]{Ga:colored}}]\label{lem-typ}
In $\Repsl$, typical, simple, and projective are equivalent for $V(k,a)$.
\end{lemma}

We have the following typical tensor product formula and their highest weight
vectors in the alternative order. Compare with \cite[Lemma~3.7]{Ga:colored}
which uses the standard order.

\begin{lemma}
\label{lemma-tensor}
The tensor product formula
\begin{equation}
\label{Vten}
V(k,a) \otimes V(0,b) \cong
V(k+1,c) \oplus V(k,c) \oplus V(k,c+1) \oplus (1-\delta_{k,0}) V(k-1,c+1)
\end{equation}
is valid when all modules in the formula are typical, with the understanding
that the last term is omitted when $k=0$. The highest weight vectors in the
alternative order are given by
\begin{equation}
\label{eq-high-v}
\begin{aligned}
V(k+1,c)&:& m_1&=v_{0,1}\otimes v_1,\\
V(k,c)&:& m_2&=v_{0,1}\otimes v_2+(-1)^kq^{-b}v_{0,2}\otimes v_1,\\
V(k,c+1)&:& m_3&=v_{0,1}\otimes v_3+(-1)^k\frac{q^{b+1}-q^{-b-1}}{q^{2(a+k
+1)}-1}v_{0,3}\otimes v_1,\\
V(k-1,c+1)&:& m_4&=v_{0,1}\otimes v_4+(-1)^kq^{-b-1} v_{0,2}\otimes v_3
+(-1)^{k+1}q^{b-1}\delta v_{0,3}\otimes v_2\\
&&&\qquad+q^{-1}\delta v_{0,4}\otimes v_1
+\frac{\CC{1}}{\CC{k}}(q^k\delta+q^{-k})v_{1,1}\otimes v_1,\\
&&&\qquad\text{where }\delta=\frac{1-q^{-2b}}{q^{2(a+k+1)}-1}.
\end{aligned}
\end{equation}
\end{lemma}

\begin{proof}
The irreducible factors can be determined by characters. See
\cite[Lemma~1.3]{GP:multi}. Explicit calculation shows that the vectors $m_i$
are annihilated by $F_2$ and $E_{12}$ and have the expected weights.
\end{proof}

\subsection{Basic tensors}
\label{sec-basic-maps}

Let $\cat$ denote the full subcategory of $\Repsl$ whose objects are nontrivial
tensor products of typical $V(0,a),a\in\colorL\setminus\{0,-1\}$. In this
section we define some basic tensors of the category $\cat$ which will be the
images of our $\gsl$-webs under the Reshetikhin--Turaev construction. Recall we
defined $a^\ast=-a-1$ in \eqref{eq-color-flip}. There is an isomorphism
$\phi_a:V(0,a)^\ast\to V(0,a^\ast)$ given by
\begin{equation}
\label{eq-dual-iso}
\begin{alignedat}{2}
\phi_a(v_1^\ast)&=-q^{-a-1}v_4,&
\phi_a(v_2^\ast)&=-q^{-1}v_3,\\
\phi_a(v_3^\ast)&=-q^{-2a-3}v_2,&\qquad
\phi_a(v_4^\ast)&=q^{-a-2}v_1,
\end{alignedat}
\end{equation}
where $v_i^\ast$ are the dual basis to $v_i$. The scalars are chosen 
so that the following maps
\begin{equation}
\begin{aligned}
b_a&:\groundK\xrightarrow{\coev}V(0,a)\otimes V(0,a)^\ast
\xrightarrow{1\otimes\phi_a}V(0,a)\otimes V(0,a^\ast),\\
d_a&:V(0,a^\ast)\otimes V(0,a)\xrightarrow{\phi_a^{-1}\otimes1}
V(0,a)^\ast\otimes V(0,a)\xrightarrow{\ev}\groundK
\end{aligned}
\end{equation}
agree with those in \cite{GP:multi} (written with a bar), where the 
evaluation
and coevaluation maps are standard
\begin{equation}
\ev(f\otimes v)=f(v),\quad f\in V^\ast,v\in V,\qquad
\coev(1)=\sum_i v_i\otimes v_i^\ast.
\end{equation}
Explicitly, the maps $b_a$ and $d_a$ are given by
\begin{gather}
\label{eq-ba}
b_a(1)=-q^{-a-1}v_1\otimes v_4-q^{-1}v_2\otimes v_3-q^{-2a-3}v_3\otimes v_2
+q^{-a-2}v_4\otimes v_1.\\
\begin{alignedat}{2}
d_a(v_1\otimes v_4)&=q^{a+2},&
d_a(v_2\otimes v_3)&=-q^{2a+3},\\
d_a(v_3\otimes v_2)&=-q,&\qquad
d_a(v_4\otimes v_1)&=-q^{a+1},
\end{alignedat}
\end{gather}
and all other basis pairings are zero for $d_a$.

We recall two more maps from \cite{GP:multi} except $\gamma_-$ is 
renormalized so that the skein relations are as simple as possible.
\begin{equation}\label{eq-gamma}
\gamma^{a,b}_+:V(0,a+b)\to V(0,a)\otimes V(0,b),\qquad
\gamma^{a,b}_-:V(0,a+b+1)\to V(0,a)\otimes V(0,b)
\end{equation}
given by
\begin{align}
&\begin{aligned}
\gamma^{a,b}_+(v_1)&=q^b v_1\otimes v_2+v_2\otimes v_1,\\
\gamma^{a,b}_+(v_2)&=v_2\otimes v_2,\\
\gamma^{a,b}_+(v_3)&=q^{b+1}v_1\otimes v_4+v_2\otimes v_3+q^{2b}v_3\otimes 
v_2-q^b v_4\otimes v_1,\\
\gamma^{a,b}_+(v_4)&=v_2\otimes v_4+q^b v_4\otimes v_2.
\end{aligned}\\
&\begin{aligned}
\gamma^{a,b}_-(v_1)
&=-q^{-1}\CC{b+1}^{-1}v_1\otimes v_3-q^{-a-2}\CC{a+1}^{-1}v_3\otimes v_1,\\
\gamma^{a,b}_-(v_2)&=
\CC{a+b+1}^{-1}\big(q^{-a-1}v_1\otimes v_4-q^{-a-2}v_4\otimes v_1\big)\\
&\quad-\CC{a+b+1}^{-1}\Big(q^{-1}\frac{\CC{a}}{\CC{b+1}}v_2\otimes v_3
+q^{-2a-3}\frac{\CC{b}}{\CC{a+1}}v_3\otimes v_2\Big),\\
\gamma^{a,b}_-(v_3)&=-q^{-1}\frac{\CC{a+b+2}}{\CC{a+1}\CC{b+1}}v_3\otimes 
v_3,\\
\gamma^{a,b}_-(v_4)&=-q^{-a-2}\CC{a+1}^{-1}v_3\otimes v_4-q^{-1}\CC{b+1}^{-1}
v_4\otimes v_3.
\end{aligned}
\end{align}
Comparing with the highest weight vectors from \eqref{eq-high-v}, we get
\begin{equation}
\gamma^{a,b}_+(v_1)=q^bm_2,\qquad
\gamma^{a,b}_-(v_1)=-q^{-1}\CC{b+1}^{-1}m_3.
\end{equation}

\subsection{Multiplicities}
\label{sub.mult}

In this section we give formulas for the multiplicities
of the decomposition of tensor products of $V(0,a)$.

We say that a list $a_1,\dotsc,a_n$ of colors is generic if the sum of any
consecutive subset of length $i$ is not any of $0,-1,\dotsc,-i$.
Anghel~\cite{Anghel} studied the decomposition of the tensor product
$V(0,a)^{\otimes n}$, which is easily generalized to a generic list
$a_1,\dotsc,a_n$
\begin{equation}
\label{otimesn}
\bigotimes_{i=1}^n V(0, a_i) \cong
\bigoplus_{\substack{k,\ell \geq 0\\k+\ell \leq n-1}}
m^{(n)}_{k,\ell} V(k, b + \ell),
\qquad \text{where }b=a_1 + \dotsb + a_n.
\end{equation}
The simple components are labeled by certain paths in a grid. Define a grid
path of length $n-1$ as a sequence
$[M_{s_1},\dotsc,M_{s_{n-1}}]:=\big(\sum_{j=1}^t M_{s_j}\big)_{t=0}^{n-1}$ 
where
\begin{equation}\label{eq-moves}
M_1=(1,0),\quad M_2=(0,0),\quad M_3=(0,1),\quad M_4=(-1,1)
\end{equation}
such that the end of each step, i.e.\ the partial sum $\sum_{j=1}^t M_{s_i}$,
stay in the first quadrant $\nats^2$. Let $\paths^{(n)}_{k,\ell}$ be the 
set of
grid path of length $n-1$ ending at $(k,\ell)$. Then the $V(k, b + \ell)$
summands correspond to grid paths in $\paths^{(n)}_{k,\ell}$. Note the 
moves are
reordered compared to \cite{Anghel} for the same reason in
Section~\ref{sec-roots}.

\begin{remark}
The grid paths here are related to the paths in the dominant chamber of the
weight lattice (see e.g.\ \cite{KK}) by a uniform shift in the second component.
Each move is the shifted weight of a basis vector. Compare
Figure~\ref{fig-basis-weights}.
\end{remark}

\begin{remark}\label{rem.proj}
\eqref{otimesn} shows that all typical $V(k,b)$ show up as a direct summand of
tensor products of $V(0,a)$. More generally, by \cite{GQS:rep}, all projective
objects of $\Repsl$ are direct summands of tensor products of $V(0,a)$. In other
words, the ideal $\mathcal{I}$ of projective objects in $\Repsl$ is generated by
typical $V(0,a)$. Therefore, it suffices to work with $\cat$ instead of
$\mathcal{I}$ in Theorem~\ref{thm.2}.
\end{remark}

The multiplicities $m^{(n)}_{k,\ell}$ can be computed recursively by
\begin{equation}
\label{mrec}
m^{(n)}_{k,\ell} = m^{(n-1)}_{k,\ell} + m^{(n-1)}_{k,\ell-1} + m^{(n-1)}_
{k-1,\ell}
+ m^{(n-1)}_{k+1,\ell-1},\qquad
m^{(1)}_{0,0}=1.
\end{equation}
Here, the multiplicities are extended by zero, that is, $m^{(n)}_{k,\ell}=0$ if
$k<0$, $\ell<0$, or $k+l\ge n$.

\begin{proposition}
\label{prop.mult}  
The multiplicities are given by
\begin{equation}
\label{eq-mult}
m^{(n)}_{k,\ell}=\frac{k+1}{n}\binom{n}{\ell}\binom{n}{k+\ell+1} \,.
\end{equation}
The formula works for $k\ge-1$ and all $\ell$ with the convention
$\binom{n}{a}=0$ if $a$ or $n-a$ is negative.
\end{proposition}

\begin{proof}
Straightforward verification by induction on $n$.
\end{proof}

\begin{corollary}
\label{cor.cat4}
For all $n \geq 1$, we have
\begin{equation}
\dim_{\groundK}(\End(V(0,a)^{\otimes n})) = m^{(2n-1)}_{0,n-1} \,.
\end{equation}
\end{corollary}

\begin{proof}
This follows by comparing Anghel's formula \cite[Thm.~1]{Anghel} 
with \eqref{eq-mult}. We also give an alternative proof. If we rewrite
\begin{equation}
E_n:=\End(V(0,a)^{\otimes n})
\cong\Hom(V(0,a),V(0,a)^{\otimes n}\otimes V(0,
a^\ast)^{\otimes(n-1)}),
\end{equation}
the tensor colors are $C=(a,\dotsc,a,a^\ast,\dotsc,a^\ast)$, which is not
generic. However, note that the tensor product formula \eqref{Vten} only depends
on the sum of colors. This means we can choose a generic list
$C'=(a'_1,\dotsc,a'_{2n-1})$ with the same sum as $C$, and
$E_n\cong\Hom(V(0,a),\bigotimes_{i=1}^{2n-1}V(0,a'_i))$ as vector spaces. The
latter has dimension $m^{(2n-1)}_{0,n-1}$ by the decomposition~\eqref{otimesn}.
\end{proof}

\subsection{Highest weight vectors}
\label{sec-hwv}

We determine the highest weight vectors for the simple components of the
tensor product \eqref{otimesn}.

Essentially by construction, each grid path $P=[M_{s_1},\dotsc,M_{s_{n-1}}]$
ending at $(k,\ell)$ determines an inclusion map
\begin{equation}
i_P:V_P\embed \bigotimes_{i=1}^n V(0, a_i),
\end{equation}
where $V_P:=V(k,b+\ell)$ for convenience. This is defined inductively using 
the
highest weight vectors from Lemma~\ref{lemma-tensor}. More precisely, let 
$i_s$
be the inclusion of the $s$-th term on the right-hand side of \eqref{Vten} by
$v_{0,1}\mapsto m_s$ as in \eqref{eq-high-v}. Define $i_P=\id$ if $P=[\,]$ 
is trivial, or else writing $P=[P',M_s]$,
\begin{equation}
\label{eq-path-incl}
i_P:V_P\xhookrightarrow{i_s}
V_{P'}\otimes V(0,a_n)
\xhookrightarrow{i_{P'}\otimes\id}
\bigotimes_{i=1}^n V(0, a_i).
\end{equation}
A quick look at the term orders of the highest weight vectors in
\eqref{eq-high-v} proves the following lemma. This is analogous to 
\cite[Thm.~2]{KK}.

\begin{lemma}\label{lemma-highest}
In the alternative order, the highest term of the image $i_P(v_{0,1})$ of the
highest weight vector for a grid path $P=[M_{s_1},\dotsc,M_{s_{n-1}}]$ is
\begin{equation}
i_P(v_{0,1})=v_1\otimes v_{s_1}\otimes\dotsm\otimes v_{s_{n-1}}+\text{lower 
terms}.
\end{equation}
\end{lemma}

\subsection{Generating functions}
\label{sec-gen-fn}

Define generating functions for the multiplicities $m^{(n)}_{k,\ell}$.
\begin{equation}
g_k(x,y)=\sum_{n=k+1}^\infty \sum_{\ell=0}^{n-1-k} m^{(n)}_{k,\ell} x^n y^\ell
\in \ints[\![x,y]\!],
\quad
G(x,y,z)=\sum_{k=0}^\infty g_k(x,y)z^k \in \ints[\![x,y,z]\!] \,.
\end{equation}

\begin{theorem}\label{thm.gen}
The total generating function $G(x,y,z)$ is given by
\begin{equation}\label{eq-gen-fn}
G(x,y,z)=\frac{xyg_0(x,y)-xz}{xz^2+(xy+x-1)z+xy},
\end{equation}
where $g_0(x,y)=G(x,y,0)$ is given by
\begin{equation}\label{eq-g0}
g_0(x,y)=\frac{2x}{1-x-xy+\sqrt{(1-x-xy)^2-4x^2y}}.
\end{equation}
\end{theorem}

The proof relies on two new recursive formulas which will be given later.

\begin{proof}
A standard calculation using Lemma~\ref{lemma-recurs} below shows that $g_0$
satisfies the quadratic equation
\begin{equation}
g_0(x,y)-x=(x+xy)g_0(x,y)+xyg_0(x,y)^2.
\end{equation}
The solution which is finite at $(x,y)=(0,0)$ is given by \eqref{eq-g0}. Then
$G(x,y,z)$ is determined by Lemma~\ref{lem.gen} below by the usual methods.
\end{proof}

The first recursive formula used in the proof is for $m^{(n)}_{0,\ell}$ that
does not use $k>0$. (In fact, this is observed in the tree count version
Lemma~\ref{lemma-tree-count} first.)

\begin{lemma}
\label{lemma-recurs}
$m^{(n)}_{0,\ell}$ satisfies the recursive relation
\begin{equation}
m^{(n)}_{0,\ell}=m^{(n-1)}_{0,\ell}+m^{(n-1)}_{0,\ell-1}
+\sum_{i=1}^\ell\sum_{j=i+1}^{n-1} m^{(j-1)}_{0,i-1}m^{(n-j)}_{0,\ell-i},
\qquad n>1.
\end{equation}
\end{lemma}

\begin{proof}
Note the extension by zero is compatible with the path count. The first move
cannot be $M_4$ by the first quadrant condition. If the first move is $M_2$,
then we still have $m^{(n-1)}_{0,\ell}$ choices for the remaining $n-1$ moves.
Similarly, $M_3$ can be followed by $m^{(n-1)}_{0,\ell-1}$ choices. Finally, if
the first move is $M_1$, let $(0,i)$ be the first time the path returns to
$k=0$. This means after the first move, the path stays in $k\ge1$ until an $M_4$
move $(1,i-1)\to(0,i)$. There are $m^{(j-1)}_{0,i-1}m^{(n-j)}_{0,\ell-i}$
choices if the $M_4$ move above is the $j$-th move.
\end{proof}

The second recursive formula is a simple rewrite of \eqref{mrec}.

\begin{lemma}
\label{lem.gen}
$g_k$ satisfies the recursive relations
\begin{align}
\label{eq-g1}
xyg_1(x,y)&=(1-x-xy)g_0(x,y)-x,\\
\label{eq-gk}
xyg_k(x,y)&=(1-x-xy)g_{k-1}(x,y)-xg_{k-2}(x,y),\qquad k\ge2.
\end{align}
\end{lemma}

\begin{proof}
Rewrite \eqref{mrec} by grouping the terms according to $k$ and shifting indices
as follows.
\begin{equation}
m^{(n)}_{k,\ell} = m^{(n+1)}_{k-1,\ell+1} - m^{(n)}_{k-1,\ell+1}
- m^{(n)}_{k-1,\ell} - m^{(n)}_{k-2,\ell+1}.
\end{equation}
The usual manipulation of generating series, using \eqref{eq-mult} for boundary
terms, proves the lemma.
\end{proof}

The following is an interesting combinatorial fact which was observed first and
motivated the proof of Lemma~\ref{lemma-recurs}.

\begin{corollary}
\label{cor.cat}
The sum
\begin{equation}
\sum_{\ell=0}^{n-1} m^{(n)}_{0,\ell} = C_n := \frac{1}{n+1}\frac{(2n)!}{n!^2}
\end{equation}
is given by the Catalan number.
\end{corollary}

\begin{proof}
The generating function of the sum is given by
\begin{equation}
g_0(x,1)=\frac{2x}{1-2x+\sqrt{(1-2x)^2-4x^2}}=\frac{1-\sqrt{1-4x}}{2x}-1.
\end{equation}
This is the well-known generating function of the Catalan numbers with the 0-th
term removed.
\end{proof}

%%%%%%%%%%%%%%%%%%%%%%%%%%%%%%%%%%%%%%%%%%%%%%%%%%%%%%%%%%%%%%%%%%%%%%%%%%%%
%%%%%%%%%%%%%%%%%%%%%%%%%%%%%%%%%%%%%%%%%%%%%%%%%%%%%%%%%%%%%%%%%%%%%%%%%%%%

\section{Skein theory}

\subsection{Skein relations and the skein module}
\label{sec-skein}

We now give the precise definition of $\gsl$-webs following and generalizing
\cite{GP:multi}. Recall the definition of $\gsl$-webs from the Introduction.

Suppose $M$ is an oriented 3-manifold, and $V\subset\partial M$ is a set of
colored points on the boundary. Given a color lattice $\colorL$, let 
$\webcat(M,V)$ denote the
free $\groundR$-module generated by the set of isotopy classes of
$\gsl$-webs in $(M,V)$. Before we define the 
$\colorL$-colored $\gsl$-skein module $\skein(M,V)$, we first
define a perturbed skein module. Consider the perturbed color lattice 
and ground ring
\begin{equation}
\label{lambdae}
\colorL_\varepsilon=\colorL\oplus\ints\varepsilon,\qquad
R_\varepsilon:=R_{\colorL_\varepsilon}
=\groundR[q^{\pm\varepsilon},q^{\pm2a_i\varepsilon},q^{\pm2\varepsilon^2}]
[\CC{b}^{-1},b\in\colorL_\varepsilon\setminus\colorL],
\end{equation}
where $1,a_1,\dotsc,a_r$ is a basis of the color lattice $\colorL$. By
substitution $\varepsilon=0$, we mean the projection
$\colorL_\varepsilon\onto\colorL$ and the partially defined map
$R_\varepsilon\dashrightarrow\groundR$ setting the new generators of
$R_\varepsilon$ to $1$.

Choose perturbed colors $V_\varepsilon$ that recover $V$ at $\varepsilon=0$. Let
$\webcat_\varepsilon(M,V_\varepsilon)$ denote the free $R_\varepsilon$-module
generated by the set of isotopy classes of $\gsl$-webs in $(M,V_\varepsilon)$.
Define the perturbed skein module $\skein_\varepsilon(M,V_\varepsilon)$ as the
quotient of $\webcat_\varepsilon(M,V_\varepsilon)$ by the skein
relations~\eqref{eq-split-link}, \eqref{eq-bubble}, \eqref{eq-IH},
\eqref{eq-crossing}. As usual, in the relations, the webs differ only inside a
ball, and the projections of this ball is shown in the equations.

\begin{remark}\label{rem-rel-conds}
Note each relation has the condition that all webs have typical colors. For
example, $a=b$ is not allowed in the crossing relation \eqref{eq-crossing} due
to the first diagram. (Additionally, $a-b\ne-1$ and $a+b\ne0,-1,-2$.) The $a=b$
case is allowed in the real $\gsl$-skein module but takes a different form
\eqref{eq-equal-cross}.
\end{remark}

For the ball $M=B^3$, let $I_\varepsilon(B^3,V_\varepsilon)$ denote the
submodule generated by the skein relations above, so by definition,
$\skein_\varepsilon(B^3,V_\varepsilon)=\webcat_\varepsilon(B^3,V_\varepsilon)/I_\varepsilon(B^3,V_\varepsilon)$.
As before, there is a maximal $\groundR$-submodule
$I'(B^3,V_\varepsilon)\subset I_\varepsilon(B^3,V_\varepsilon)$ where the
substitution $\varepsilon=0$ is well-defined, and let
$I(B^3,V)\subset\webcat(B^3,V)$ denote the sum of the substitutions $\varepsilon=0$
for all choices of perturbations $V_\varepsilon$ of $V$.

A $\colorL$-colored $\gsl$-skein relation in $(M,V)$ is a relation $\sum_i
c_i\alpha_i=0$ where $c_i\in \groundR$, and $\alpha_i\in\webcat(M,V)$ are
identical outside a ball $B^3\subset M$, such that
\begin{equation}
\sum_i c_i(\alpha_i\cap B^3)\in I(B^3,V')
\end{equation}
where $V'\subset\partial B^3$ is the set of colored points common to all terms.
This is the precise formulation of the perturbation in
Section~\ref{sec-intro-skein}. Note the perturbed webs $D_i(\varepsilon)$ in the
ball need not extend to perturbations of the whole $\alpha_i$.

Finally, the $\colorL$-colored $\gsl$-skein module of $(M,V)$, denoted by
$\skein(M,V)$, is the quotient of $\webcat(M,V)$ by all $\colorL$-colored
$\gsl$-skein relations. By choosing the trivial perturbation $V_\varepsilon=V$,
the perturbed skein relations \eqref{eq-split-link}--\eqref{eq-crossing} are
also $\gsl$-skein relations in $(M,V)$.

\begin{remark}
We point out some differences between our skein relations and those in
\cite{GP:multi}.
The vertices here are renormalized versions according to
\begin{equation}
\begin{linkdiag}
\fill[gray!20] (-0.3,-0.3)rectangle(1.3,1.3);
\webnode{0.5,0.5}{V}{+};
\draw[thick] (V) -- (-0.2,1.3) (V) -- (1.2,1.3) (V) -- (0.5,-0.3);
\end{linkdiag}
=
\begin{linkdiag}
\fill[gray!20] (-0.3,-0.3)rectangle(1.3,1.3);
\begin{scope}[thick]
\draw (0.5,0.5)node[right]{\stsize$+$} -- (-0.2,1.3);
\draw (0.5,0.5) -- (1.2,1.3);
\draw (0.5,0.5) -- (0.5,-0.3);
\end{scope}
\end{linkdiag},
\qquad
\begin{linkdiag}
\fill[gray!20] (-0.3,-0.3)rectangle(1.3,1.3);
\webnode{0.5,0.5}{V}{-};
\begin{scope}[thick]
\draw[-o-={0.8}{>}] (V) --node[pos=0.4,left]{$a$} (-0.2,1.3);
\draw[-o-={0.8}{>}] (V) --node[pos=0.4,right]{$c$} (1.2,1.3);
\draw[-o-={0.7}{>}] (V) --node[left]{$b$} (0.5,-0.3);
\end{scope}
\end{linkdiag}
=\frac{1}{\CC{a}\CC{b}\CC{c}}
\begin{linkdiag}
\fill[gray!20] (-0.3,-0.3)rectangle(1.3,1.3);
\begin{scope}[thick]
\draw[-o-={0.8}{>}] (0.5,0.5)node[right]{\stsize$-$} -- (-0.2,1.3);
\draw[-o-={0.8}{>}] (0.5,0.5) -- (1.2,1.3);
\draw[-o-={0.8}{>}] (0.5,0.5) -- (0.5,-0.3);
\end{scope}
\end{linkdiag}.
\end{equation}
The left-hand sides are our vertices, and the right-hand sides are theirs. This
choice eliminates all coefficients in $I=H$ \eqref{eq-IH}.

The second bubble relation in \eqref{eq-bubble} is not explicitly stated in \cite{GP:multi}.

Of the relations in \cite[Prop.~4.6]{GP:multi}, the last relation does not
make sense for non-closed webs, and the rest are equivalent to or follows from
our relations, so they are not listed as defining relations here. Some of them
are derived in the lemma below.

Our perturbation is more general than theirs. Their perturbation is global, but
here as a skein relation, we allow local perturbation, meaning perturbation of
colors inside a small ball without matching the colors outside the ball. For
example, this is why we can perturb the equal color crossing without worrying
whether the two strands are parts of the same edge.
\end{remark}

\subsection{Useful relations}

\begin{lemma}\label{lemma-more-rel}
The following relations are consequences of the defining relations.
\begin{enumerate}
\item Opposite crossing:
\begin{equation}\label{eq-crossing2}
\begin{linkdiag}
\fill[gray!20] (-0.3,-0.3)rectangle(1.3,1.3);
\begin{scope}[thick]
\begin{knot}
\strand (-0.2,1.3) -- (1.2,-0.3);
\strand (1.2,1.3) -- (-0.2,-0.3);
\end{knot}
\path[-o-={0.8}{>}] (1.2,1.3) -- node[pos=0.7,left]{$a\vphantom{b}$} (-0.2,-0.3);
\path[-o-={0.8}{>}] (-0.2,1.3) -- node[pos=0.7,right]{$b$} (1.2,-0.3);
\end{scope}
\end{linkdiag}
=q^{2ab+a+b}\left(\CC{a-b}\Hskein{+}{-}-q^{-b}\CC{a}\Iskein{-}{+}
+q^{a+1}\CC{b+1}\Iskein{+}{-}\right).
\end{equation}
\item Square breaking:
\begin{equation}\label{eq-square}
\begin{linkdiag}[1]
\fill[gray!20] (0,0) rectangle (2,2);
\webnode{0.5,1.5}{A}{+};
\webnode{1.5,1.5}{B}{-};
\webnode{0.5,0.5}{C}{+};
\webnode{1.5,0.5}{D}{-};
\begin{scope}[thick,inner sep=3pt]
\draw[-o-={0.5}{<}] (A) -- node[above left,pos=0.9]{$e$} (B);
\draw[-o-={0.5}{<}] (C) -- node[below left,pos=0.9]{$-e$} (D);
\draw (0.1,2) -- (A) -- (C) -- (0.1,0);
\draw (1.9,2) -- (B) -- (D) -- (1.9,0);
\end{scope}
\end{linkdiag}
=-\CC{e}^{-2}\iddiag.
\end{equation}
\item Equal color crossing:
\begin{equation}\label{eq-equal-cross}
\begin{aligned}
\begin{linkdiag}
\fill[gray!20] (-0.3,-0.3)rectangle(1.3,1.3);
\begin{scope}[thick]
\begin{knot}
\strand (1.2,1.3) -- (-0.2,-0.3);
\strand (-0.2,1.3) -- (1.2,-0.3);
\end{knot}
\path[-o-={0.8}{>}] (1.2,1.3) -- node[pos=0.7,left]{$a$} (-0.2,-0.3);
\path[-o-={0.8}{>}] (-0.2,1.3) -- node[pos=0.7,right]{$a$} (1.2,-0.3);
\end{scope}
\end{linkdiag}
&=q^{-2a^2-2a}\left(-\iddiag-q^a\CC{a}\Iskein{-}{+}
+q^{-a-1}\CC{a+1}\Iskein{+}{-}\right).\\
\begin{linkdiag}
\fill[gray!20] (-0.3,-0.3)rectangle(1.3,1.3);
\begin{scope}[thick]
\begin{knot}
\strand (-0.2,1.3) -- (1.2,-0.3);
\strand (1.2,1.3) -- (-0.2,-0.3);
\end{knot}
\path[-o-={0.8}{>}] (1.2,1.3) -- node[pos=0.7,left]{$a$} (-0.2,-0.3);
\path[-o-={0.8}{>}] (-0.2,1.3) -- node[pos=0.7,right]{$a$} (1.2,-0.3);
\end{scope}
\end{linkdiag}
&=q^{2a^2+2a}\left(-\iddiag-q^{-a}\CC{a}\Iskein{-}{+}
+q^{a+1}\CC{a+1}\Iskein{+}{-}\right).
\end{aligned}
\end{equation}
\item Reidemeister I:
\begin{equation}\label{eq-RI}
\begin{linkdiag}
\fill[gray!20] (-0.1,-0.1)rectangle(1.1,1.1);
\begin{knot}
\strand[thick] (1,0.5) ..controls +(0,0.4) and +(0,1).. (0.3,-0.1);
\strand[thick] (0.3,1.1) ..controls +(0,-1) and +(0,-0.4).. (1,0.5);
\end{knot}
\path[thick,-o-={0.9}{>}] (1,0.5) ..controls +(0,0.4) and +(0,1)..
  node[pos=0.9,left]{$a$} (0.3,-0.1);
\end{linkdiag}
=q^{-2a^2-2a}
\begin{linkdiag}
\fill[gray!20] (-0.1,-0.1)rectangle(1.1,1.1);
\draw[thick] (0.3,1.1) -- (0.3,-0.1);
\end{linkdiag},
\qquad
\begin{linkdiag}
\fill[gray!20] (-0.1,-0.1)rectangle(1.1,1.1);
\begin{knot}
\strand[thick] (0.3,1.1) ..controls +(0,-1) and +(0,-0.4).. (1,0.5);
\strand[thick] (1,0.5) ..controls +(0,0.4) and +(0,1).. (0.3,-0.1);
\end{knot}
\path[thick,-o-={0.9}{>}] (1,0.5) ..controls +(0,0.4) and +(0,1)..
  node[pos=0.9,left]{$a$} (0.3,-0.1);
\end{linkdiag}
=q^{2a^2+2a}
\begin{linkdiag}
\fill[gray!20] (-0.1,-0.1)rectangle(1.1,1.1);
\draw[thick] (0.3,1.1) -- (0.3,-0.1);
\end{linkdiag}.
\end{equation}
\end{enumerate}
\end{lemma}

\begin{proof}
\begin{enumcc}[(1)]
\item \eqref{eq-crossing2} is obtained from \eqref{eq-crossing} by reversing the
orientation on one of the strands and applying the last $I=H$ relation in
\eqref{eq-IH}.

\item Start with the Reidemeister II move
\begin{equation}
\begin{linkdiag}
\fill[gray!20] (-0.3,-0.3)rectangle(1.3,1.3);
\begin{scope}[thick]
\draw[-o-={0.8}{>}] (-0.2,1.3)..controls(0.2,0.5)..(-0.2,-0.3) node[pos=0.8,right]{$a$\vphantom{$b$}};
\draw[-o-={0.8}{>}] (1.2,1.3)..controls(0.8,0.5)..(1.2,-0.3) node[pos=0.8,left]{$b$};
\end{scope}
\end{linkdiag}
=
\begin{linkdiag}
\fill[gray!20] (-0.3,-0.3)rectangle(1.3,1.3);
\begin{knot}
\strand[thick] (-0.2,1.3)..controls(1.3,0.5)..(-0.2,-0.3);
\strand[thick] (1.2,1.3)..controls(-0.3,0.5)..(1.2,-0.3);
\end{knot}
\end{linkdiag}
\end{equation}
and resolve the crossings on the right-hand side. After simplification using
$I=H$ \eqref{eq-IH} and bubble \eqref{eq-bubble}, the result is
\eqref{eq-square} with $e=b-a$. Note this has the condition $a-b\ne0,-1$ and
$a+b\ne0,-1,-2$ as mentioned before.

With the same conditions on $a,b$, now suppose $d:=e-(b-a)\ne0,-1$. This case
follows from the special case $e=b-a$ by creating a bubble on the $e$ edge and
transferring it to the $-e$ edge using $I=H$ \eqref{eq-IH} as follows.
\begin{equation}\label{eq-sq-bubble}
\begin{linkdiag}[1]
\fill[gray!20] (0,0) rectangle (2,2);
\webnode{0.5,1.5}{A}{+};
\webnode{1.5,1.5}{B}{-};
\webnode{0.5,0.5}{C}{+};
\webnode{1.5,0.5}{D}{-};
\begin{scope}[thick,inner sep=3pt]
\draw[-o-={0.5}{<}] (A) -- node[above left,pos=0.9]{$e$} (B);
\draw[-o-={0.5}{<}] (C) -- node[below left,pos=0.9]{$-e$} (D);
\draw (0.1,2) -- (A) -- (C) -- (0.1,0);
\draw (1.9,2) -- (B) -- (D) -- (1.9,0);
\end{scope}
\end{linkdiag}
\longrightarrow
\begin{linkdiag}[1.25]
\fill[gray!20] (0,0) rectangle (3.5,2.5);
\webnode{0.5,2}{A}{+};
\webnode{3,2}{B}{-};
\webnode{0.5,0.5}{C}{+};
\webnode{3,0.5}{D}{-};
\webnode{1.1,2}{L}{+};
\webnode{2.4,2}{R}{-};
\begin{scope}[thick,inner sep=3pt]
\draw (A) -- (L) (R) -- (B);
\draw[-o-={0.5}{<}] (C) -- node[below]{$-e$} (D);
\draw (0.1,2.5) -- (A) -- (C) -- (0.1,0);
\draw (3.4,2.5) -- (B) -- (D) -- (3.4,0);
\draw[-o-={0.5}{<}] (L) -- node[above]{$b-a$} (R);
\draw[-o-={0.5}{<}] (L) ..controls +(0.2,-0.4) and +(-0.2,-0.4).. node[below]{$d$}(R);
\end{scope}
\end{linkdiag}
\longrightarrow
\begin{linkdiag}[1.5]
\fill[gray!20] (0,0) rectangle (2.5,3);
\webnode{0.5,2.5}{A}{+};
\webnode{2,2.5}{B}{-};
\webnode{0.5,1.5}{L}{+};
\webnode{2,1.5}{R}{-};
\webnode{0.5,0.5}{C}{+};
\webnode{2,0.5}{D}{-};
\begin{scope}[thick,inner sep=3pt]
\draw[-o-={0.5}{<}] (A) -- node[above]{$b-a$} (B);
\draw[-o-={0.5}{<}] (L) -- node[below]{$d$} (R);
\draw[-o-={0.5}{<}] (C) -- node[below]{$-e$} (D);
\draw (0.1,3) -- (A) -- (L) -- (C) -- (0.1,0);
\draw (2.4,3) -- (B) -- (R) -- (D) -- (2.4,0);
\end{scope}
\end{linkdiag}
\longrightarrow
\begin{linkdiag}[1.25]
\fill[gray!20] (0,0) rectangle (3.5,2.5);
\webnode{0.5,2}{A}{+};
\webnode{3,2}{B}{-};
\webnode{0.5,0.5}{C}{+};
\webnode{3,0.5}{D}{-};
\webnode{1.1,0.5}{L}{+};
\webnode{2.4,0.5}{R}{-};
\begin{scope}[thick,inner sep=3pt]
\draw[-o-={0.5}{<}] (A) -- node[above]{$b-a$} (B);
\draw (C) -- (L) (R) -- (D);
\draw (0.1,2.5) -- (A) -- (C) -- (0.1,0);
\draw (3.4,2.5) -- (B) -- (D) -- (3.4,0);
\draw[-o-={0.5}{<}] (L) -- node[below]{$-e$} (R);
\draw[-o-={0.5}{<}] (L) ..controls +(0.2,0.4) and +(-0.2,0.4).. node[above]{$d$}(R);
\end{scope}
\end{linkdiag}
\longrightarrow
\begin{linkdiag}[1]
\fill[gray!20] (0,0) rectangle (2.5,2);
\webnode{0.5,1.5}{A}{+};
\webnode{2,1.5}{B}{-};
\webnode{0.5,0.5}{C}{+};
\webnode{2,0.5}{D}{-};
\begin{scope}[thick,inner sep=3pt]
\draw[-o-={0.5}{<}] (A) -- node[above]{$b-a$} (B);
\draw[-o-={0.5}{<}] (C) -- node[below]{$a-b$} (D);
\draw (0.1,2) -- (A) -- (C) -- (0.1,0);
\draw (2.4,2) -- (B) -- (D) -- (2.4,0);
\end{scope}
\end{linkdiag}
\end{equation}

Finally, the extra restrictions can be removed by perturbation.

\item The two relations are similar, so we only prove the first one.

Perturb the lower strand to $b=a-\varepsilon$ and apply the crossing relation
\eqref{eq-crossing}. The last two terms have finite coefficients at
$\varepsilon=0$, so we only need to focus on the first term and prove
\begin{equation}\label{eq-0break}
\left[\CC{\varepsilon}
\begin{linkdiag}
\fill[gray!20] (-0.5,-0.5)rectangle(1.5,1.5);
\webnode{0,0.5}{A}{+};
\webnode{1,0.5}{B}{-};
\begin{scope}[thick]
\draw[-o-={0.5}{>}] (A) -- node[right]{$a$} (-0.4,-0.5);
\draw (B) -- (1.4,-0.5);
\draw (A) -- (-0.4,1.5);
\draw[-o-={0.6}{<}] (B) -- node[left]{$a$} (1.4,1.5);
\draw[-o-={0.6}{<}] (A) -- node[above]{$\varepsilon$} (B);
\end{scope}
\end{linkdiag}
\right]_{\varepsilon=0}=-\iddiag.
\end{equation}
We use the same bubble trick as \eqref{eq-sq-bubble}.
\begin{equation}
\Hskein{+}{-}\longrightarrow
\begin{linkdiag}
\fill[gray!20] (-0.5,-0.5)rectangle(2.5,1.5);
\webnode{0,0.5}{A}{+};
\webnode{2,0.5}{B}{-};
\webnode{0.6,0.5}{L}{+};
\webnode{1.4,0.5}{R}{-};
\begin{scope}[thick]
\draw (-0.4,-0.5) -- (A) -- (-0.4,1.5);
\draw (2.4,-0.5) -- (B) -- (2.4,1.5);
\draw (A) -- (L) (R) -- (B);
\draw (L) ..controls +(0.2,0.4) and +(-0.2,0.4).. (R);
\draw (L) ..controls +(0.2,-0.4) and +(-0.2,-0.4).. (R);
\end{scope}
\end{linkdiag}
\longrightarrow
\begin{linkdiag}[1]
\fill[gray!20] (0,0) rectangle (2,2);
\webnode{0.5,1.5}{A}{+};
\webnode{1.5,1.5}{B}{-};
\webnode{0.5,0.5}{C}{+};
\webnode{1.5,0.5}{D}{-};
\begin{scope}[thick]
\draw (A) -- (B);
\draw (C) -- (D);
\draw (0.1,2) -- (A) -- (C) -- (0.1,0);
\draw (1.9,2) -- (B) -- (D) -- (1.9,0);
\end{scope}
\end{linkdiag}
\end{equation}
By carefully choosing the colors on the new edges of the bubble, the square has
a finite coefficient in the limit $\varepsilon=0$, and it breaks by
\eqref{eq-square}.

\item Resolve the crossing and then apply split link \eqref{eq-split-link} and
bubble \eqref{eq-bubble}.\qedhere
\end{enumcc}
\end{proof}

\begin{remark}
Skein modules are usually defined using 3d isotopy classes, so Reidemeister II
and III moves are automatically satisfied for diagrams. Alternatively, if $M$ is
a thickened surface $\Sigma\times[0,1]$, one can define a module using isotopy
classes of diagrams and impose relations. The proof of the lemma shows that
Reidemeister II is equivalent to square breaking. In addition, the proof of
Theorem~\ref{thm.1} shows that the same set of skein relations plus square
breaking on diagrams give the same skein module. In particular, Reidemeister III
and the following move
\begin{equation}
\begin{linkdiag}
\fill[gray!20] (-0.3,-0.3)rectangle(1.3,1.3);
\begin{knot}
\strand[thick] (0.5,0.5) -- (-0.2,1.3) (0.5,0.5) -- (1.2,1.3) (0.5,0.5) -- (0.5,-0.3);
\strand[thick] (-0.3,0.5)..controls (0.5,1.1)..(1.3,0.5);
\end{knot}
\end{linkdiag}
=
\begin{linkdiag}
\fill[gray!20] (-0.3,-0.3)rectangle(1.3,1.3);
\begin{knot}
\strand[thick] (0.5,0.5) -- (-0.2,1.3) (0.5,0.5) -- (1.2,1.3) (0.5,0.5) -- (0.5,-0.3);
\strand[thick] (-0.3,0.5)..controls (0.5,-0.1)..(1.3,0.5);
\end{knot}
\end{linkdiag}
\end{equation}
can be derived from the other relations. This fact can also be proved by direct
but extremely long calculation.
\end{remark}

So far, to resolve crossings, we have \eqref{eq-crossing}, \eqref{eq-crossing2},
and \eqref{eq-equal-cross}. Now we deal with the remaining cases $a-b=-1$ or
$a+b=0,-1,-2$ as mentioned in Remark~\ref{rem-rel-conds}.

\begin{lemma}\label{lem-all-cross-type}
All crossings can be resolved.
\end{lemma}

\begin{proof}
Note the remaining cases are mutually exclusive. In fact, if $a+b=0,-1,-2$, then
$a-b$ cannot be $1,0,-1$. Otherwise, one of $a,b$ would be atypical or $\pm1/2$,
which are forbidden.
\begin{enumcc}[(1)]
\item If $a-b=-1$, by flipping the orientations of both strands, we get
$a^\ast-b^\ast=1$, which can be reduced.
\item If $a+b=0$, flip the $a$ strand to get a crossing with $a'=b$,
$b'=a^\ast$. Then $a'-b'=a+b+1=1$ and $a'+b'=b-a-1\ne0,-1,-2$, so the new
crossing can be reduced.
\item If $a+b=-1$, i.e.\ $b=a^\ast$, flip the $a$ strand to an equal color
crossing, which can be reduced.
\item If $a+b=-2$, flip the $b$ strand to get a crossing with $a'=b^\ast$,
$b'=a$. Then $a'-b'=-a-b-1=1$ and $a'+b'=a-b-1\ne0,-1,-2$, so the new crossing
can be reduced.\qedhere
\end{enumcc}
\end{proof}

\subsection{Categories}

We now introduce three auxiliary ribbon categories that will be used 
in the proof of Theorem~\ref{thm.1}.

The first is the category $\webcat$ of $\gsl$-webs in $B^3=[0,1]^3$ whose
objects are finite sequences of typical colors, and the $\Hom$ spaces are free
$\groundR$-modules generated by $\gsl$-webs with the following boundary points.
For each $n>0$, fix a choice of $n$ distinct points
$O_n\subset\{\frac{1}{2}\}\times[0,1]$ inside the square $[0,1]^2$. Let
$V_{m,n}=O_m\times\{0\}\sqcup O_n\times\{1\}\subset B^3$. Then
$\Hom_\webcat(\vec{b},\vec{a})$ is defined as $\webcat(B^3,V_{m,n})$ where
$\vec{b}=(b_1,\dotsc,b_m)$ and $\vec{a}=(a_1,\dotsc,a_n)$ are objects of
$\webcat$, and the colors of $V_{m,n}$ are $\vec{a}$ on the top and
$(b_1^\ast,\dotsc,b_m^\ast)$ on the bottom. See Figure~\ref{fig-skcat-color} for
an example. By convention, the diagrams are front projections, i.e.\ onto the
last two components. Composition is given by vertical stacking. Moreover,
$\webcat$ is a ribbon category, where the monoidal structure is horizontal
juxtaposition, and the dualities are caps and cups.

\begin{figure}
\centering
\begin{tikzpicture}
\fill[gray!20] (0,0) rectangle (3,2);
\draw[thick] (0,0) -- (3,0) (0,2) -- (3,2);
\foreach \x/\l in {0.6/$a_1$, 1.5/$a_2$, 2.4/$a_3$} {
\fill (\x,2) circle[radius=0.1] node[above]{\l};
\draw[thick,-o-={0.6}{>}] (\x,2) -- node[left]{\l} (\x,1.2);
}
\foreach \x/\l in {1/b_1, 2/b_2} {
\fill (\x,0) circle[radius=0.1] node[below]{$\l^\ast$};
\draw[thick,-o-={0.6}{>}] (\x,0.8) -- node[left]{$\l$} (\x,0);
}
\draw[fill=white] (0.2,0.8) rectangle (2.8,1.2);
\end{tikzpicture}
\caption{Webs in $\Hom_\webcat((b_1,b_2),(a_1,a_2,a_3))$.}
\label{fig-skcat-color}
\end{figure}

The next ribbon category $\skeinballcat$ has the same objects as 
$\webcat$ and morphisms the $\gsl$-skein module of $B^3$.

Recall that $\cat$ denotes the full subcategory of $\Repsl$ whose objects are
nontrivial tensor products of typical $V(0,a)$. Using the isomorphism $\phi_a$
from Section~\ref{sec-basic-maps}, it follows that $\cat$ inherits a ribbon
structure from $\Repsl$. The last ribbon category $\ribcat_\cat$ is the category
of $\cat$-colored ribbon graphs, equipped with a functor\footnote{Technically,
we need to impose cyclic symmetry of the coupons and the color equivalence
\eqref{eq-color-flip} on $\ribcat_\cat$, but these hold in all our categories
and does not affect our statements or proofs, as explained in
\cite[Sec.~4]{GP:multi}.}
\begin{equation}
\label{w2r}  
\webcat \to \ribcat_\cat
\end{equation}
as well as the Reshetikhin--Turaev functor~\cite{RT} 
\begin{equation}
\label{r2c}  
\RT: \ribcat_\cat \to \cat.
\end{equation}

The functor of Equation~\eqref{w2r} is defined as follows. Given an $\gsl$-web,
choose an arbitrary orientation on each edge. Thicken the new graph into a
ribbon graph and color the edges in the obvious way. The 3-valent vertices are
coupons colored by maps $\gamma^{a,b}_\pm$ from \eqref{eq-gamma}, composed with
the isomorphism $\phi_a$ from \eqref{eq-dual-iso} if necessary.

The graphical convention is the following. The domain of a morphism is the
bottom, and the codomain is the top. A downward orientation is colored by the
original module, and an upward one is colored by the dual. For example, the
ribbon graph with a single coupon $\gamma^{a,b}_\pm$ looks like
$\begin{linkdiag}[0.6] \fill[gray!20] (0,0)rectangle(1.2,1.2);
\webnode{0.6,0.6}{A}{\pm};
\begin{scope}[thick]
\draw[-o-={0.6}{>}] (0.1,1.2) -- (A.120);
\draw[-o-={0.6}{>}] (1.1,1.2) -- (A.60);
\draw[-o-={0.7}{>}] (A) -- (0.6,0);
\end{scope}
\end{linkdiag}$.

Essentially by construction, we have a functor $F:\skeinballcat\to\cat$
completing the commutative diagram
\begin{equation}
\label{cd}
\begin{tikzcd}
\webcat \arrow[r] \arrow[d] & 
\skeinballcat \arrow[d,dashed,"F"] \\
\ribcat_\cat \arrow[r,"\RT"] & \cat
\end{tikzcd}
\end{equation}

Theorem~\ref{thm.1} implies the following.

\begin{theorem}
\label{thm.3}
We have an equivalence of categories
\begin{equation}
\label{b2c}
F: \skeinballcat \to \cat
\end{equation}
\end{theorem}

\begin{proof}
$F$ is a bijection on objects by construction and a bijection on morphisms by
\eqref{eq-Fbar} and dualities.
\end{proof}

\begin{proof}(of Theorem~\ref{thm.2} from Theorem~\ref{thm.3})
We have:
\begin{equation}
\skein(M,V)\otimes_{\groundR}\groundK
\simeq \Sk_{\skeinballcat}(M,V) \otimes_{\groundR}\groundK
\simeq \Sk_\cat(M,V) \simeq \skein_\mathcal{I}(M,V)
\end{equation}  
where the first isomorphism is tautological, the second one follows from
Theorem~\ref{thm.3}, and the last one is explained in Remark~\ref{rem.proj}.
\end{proof}

\subsection{Outline of the proof of Theorem~\ref{thm.1}}

By duality, we can consider
\begin{equation}\label{eq-F1n}
F_{1,n}:\skein(B^3,V_{1,n})\otimes_{\groundR}\groundK\to
\Hom_{\Uqsl}\Big(V(0,b),\bigotimes_{i=1}^n V(0,a_i)\Big)
\end{equation}
instead, where the colors of $V_{1,n}$ are $b^\ast$ and $a_1,\dotsc,a_n$.

Recall that the codomain in \eqref{eq-F1n} has a $\groundK$-basis $i_P$ labeled
by grid paths $P\in\paths^{(n)}_{0,\ell}$ defined in Section~\ref{sec-hwv}.
Therefore, our goal is to find a similar basis (over $\groundR$) for the skein
module $\skein(B^3,V_{1,n})$.

A common method of finding a basis for a skein module is to use confluence
theory \cite{SW}, where a reduction rule with nice properties implies that the
irreducible diagrams form a basis. However, it is difficult to apply confluence
to $\gsl$-webs for a few reasons. The typical conditions on the skein relations
drastically increase the cases of overlapping relations to check. The
perturbation relation can also have arbitrary complex local picture.

Nevertheless, we describe a reduction algorithm in Section~\ref{sec-reduce} for
the skein module that suggests a basis consisting of trees $\{T_P\mid
P\subset\paths^{(n)}_{0,\ell}\}$. However, for certain colors, the trees cannot
form a basis.

\begin{example}
A simple counterexample is given by $b=a$ and $(a_i)=(a,a^\ast,a)$. By drawing
all possible trees, we can check that none of them can be typically colored.
\end{example}

The algorithm perturbs the colors so that the conditions in every step are
satisfied. We will show that a good perturbation always exists. If we can always
set $\varepsilon=0$ at the end, then we prove that the trees $\{T_P\}$ span over
$\groundR$, and they become a basis by dimension count. As explained above, this
does not always happen.

In Section~\ref{sec-proof} where we give the final part of the proof, we
conclude that $\{T_P\}$ is a basis at least generically, and use it to construct
a possibly non-tree basis for the remaining cases using an isomorphism similar
to the proof of Corollary~\ref{cor.cat4}.

\subsection{Reduction algorithm}
\label{sec-reduce}

The first half of the reduction is almost the same as
\cite[Prop.~4.6]{GP:multi}. The first step is resolving all crossings, which is
always possible by Lemma~\ref{lem-all-cross-type}. Now we assume the web is
planar. For simplicity, we also assume the web is connected.

\begin{figure}
\begin{tikzpicture}
\fill[gray!20] (-0.2,0) rectangle (2,2);
\webnode[5pt]{0.5,1}{A}{};
\node[draw,rectangle,fill=white] (T) at (1.5,1) {$\beta$};
\begin{scope}[thick]
\draw[-o-={0.5}{>}] (A) -- node[above]{$k$} (-0.2,1);
\draw[-o-={0.5}{>}] (A) to[out=90,in=90] node[above]{$a$} (T);
\draw[-o-={0.5}{<}] (A) to[out=-90,in=-90] node[below]{$b$} (T);
\end{scope}
\end{tikzpicture}
\caption{Configuration that is 1-particle reducible}\label{fig-1PR}
\end{figure}
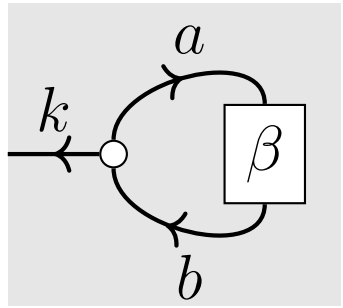

Recall that a connected graph is 1-particle irreducible if removing any internal
edge still leaves the graph connected. If the underlying graph of a planar web
is 1-particle irreducible, then it cannot contain the piece in
Figure~\ref{fig-1PR}. It is easy to see that the bubble relation
\eqref{eq-bubble} and $I=H$ relations \eqref{eq-IH} preserve 1-particle
irreducibility.

1-particle reducible webs are somewhat tricky. They are not explicitly
considered in \cite{GP:multi} since it is clear that the Reshetikhin--Turaev
invariants vanish, but it is not obvious that this follows from skein
relations. We deal with 1-particle irreducible webs first, and deal with the
reducible ones later.

We point out that without perturbation, the reduction algorithm cannot be
completed in general due to the implicit conditions of the $I=H$ relations. In
the description below, a good perturbation is implied throughout. The existence
of a good perturbation is discussed after the algorithm.

If the diagram contains polygons, choose an order of reduction. Note a polygon
has at least two sides by the typical color condition.
\begin{enumerate}
\item If the polygon has two sides, then the bubble relation \eqref{eq-bubble}
removes it.
\item If the signs around the polygon are not alternating, use the $I=H$
relations \eqref{eq-IH} of the same signs (the first two) to reduce the number
of sides.
\item Otherwise, the signs are alternating. Apply the last $I=H$ relation
\eqref{eq-IH} to produce two polygons with one fewer sides and a third diagram
with non-alternating signs for which (2) applies.
\item Repeat (2) and (3) until the current polygon is removed by (1), and then
repeat the overall process until all polygons are removed.
\end{enumerate}

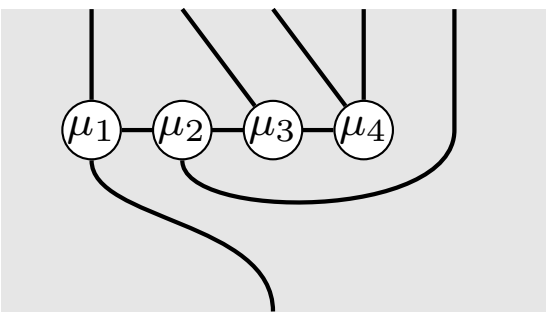
\begin{figure}
\centering
\begin{tikzpicture}[inner sep=2pt]
\fill[gray!20] (-0.6,0) rectangle (3,2);
\webnode[11pt]{0,1.2}{M1}{\mu_1};
\webnode[11pt]{1.8,1.2}{M4}{\mu_4};
\draw[thick] (0,2) -- (M1) -- (M4) -- (1.2,2) 
  (M1)..controls+(0,-0.6) and +(0,0.6)..(1.2,0) (M4) -- (1.8,2);
\webnode[11pt]{0.6,1.2}{M2}{\mu_2};
\webnode[11pt]{1.2,1.2}{M3}{\mu_3};
\draw[thick] (M2)..controls+(0,-0.6) and +(0,-0.6)..(2.4,1.2) -- (2.4,2)
  (M3) -- (0.6,2);
\end{tikzpicture}
\caption{A tree with sign list $(\mu_1,\mu_2,\mu_3),(\mu_4),\emptyset,\emptyset$.}
\label{fig-tree-signs}
\end{figure}

At this point, the diagrams are all planar trees.
We can encode the trees combinatorially as follows.  Given a tree, look at the
half-polygons, which are regions bounded by the tree and segments of the top
boundary. Order the half-polygons from left to right, and list the signs of the
vertices of the half-polygons in the counterclockwise order but excluding the
ones belonging to previous half-polygons. The result is a sequence
$S_1,\dotsc,S_{n-1}$ where each $S_i$ is a possibly empty sequence of signs
$\pm$. We call it the sign list of the tree. For example, the tree in
Figure~\ref{fig-tree-signs} has sign list
$(\mu_1,\mu_2,\mu_3),(\mu_4),\emptyset,\emptyset$. For the opposite direction,
the underlying trivalent tree can be reconstructed by $\abs{S_i}$, the vertex
colors are determined in the obvious way, and the edge colors are determined by
the boundary colors and conservation. Note this is a one-sided
inverse. Not all sequences are sign lists of trees.

An entry $S_i$ (or the corresponding half-polygon) is reduced if it is one of
the 4 possibilities: $+-$, $+$, $-$, or $\emptyset$. We reduce the trees and the
sign lists using the following procedure.
\begin{enumerate}
\item If all entries are reduced, we are done. Otherwise, find the first
non-reduced entry.
\item If the signs are not alternating, use the $I=H$ relations of the same
signs (the first two in \eqref{eq-IH}) repeatedly until the signs are
alternating.
\item If there is a pair of adjacent signs $-+$, apply the last $I=H$ relation.
As before, this produces two half-polygons with one fewer sides and a third
diagram with $-+$ replaced by $+-$.
\item Repeat from (1).
\end{enumerate}
We emphasize that if a vertex is part of a previously reduced half-polygon, it
is excluded from subsequent sign list, hence unchanged by subsequent reduction.
Note in this procedure, the choices of orders and the choice of $-+$ instead of
$+-$ are not canonical. Other choices can also work.

To see that the procedure terminates for any diagram, we just need to show that
each entry eventually becomes reduced. In the procedure, the only term that does
not have fewer sides is the last term in (3). If there are only two signs left,
then it is reduced by definition. Otherwise, the term has non-alternating signs,
which eventually reduces to fewer sides by (2). Since a sufficiently short entry
is reduced, the procedure terminates.

Now we address the perturbation.

\begin{lemma}
Suppose $\alpha$ is a planar 1-particle irreducible $\gsl$-web in $(B^3,V)$.
Then there exists a perturbation $\alpha(\varepsilon)$ such that all steps in
the algorithm above are valid.
\end{lemma}

\begin{proof}
The perturbations of $\alpha$ are parametrized by a relative $\ints$-homology
class $h\in H_1(\alpha,V)$ by the $\varepsilon$-components due to the color
equivalence \eqref{eq-color-flip} and the conservation rules
\eqref{eq-vertex-sum}. See e.g.~\cite[Prop.~4.5]{GP:multi}.

Let $\beta$ be a web used in the algorithm (including $\alpha$ itself) whose
colors depend on the initial perturbation. The $\varepsilon$-component of the
color on an edge $e$ in $\beta$ is a linear function
$\delta_\beta(e)\in\Hom(H_1(\beta,V),\ints) =H^1(\beta,V)$. Since $\alpha$ is
1-particle irreducible, so is $\beta$, and it is easy to see that
$\delta_\beta(e)$ is nontrivial.

By induction, it is easy to check that any diagram $\beta$ produced by the
algorithm is homotopy equivalent to $\alpha$ with some 2-cells attached. Let
$i_\beta:\alpha\to\beta$ be the restriction of the homotopy equivalence.

Now suppose $h\in H_1(\alpha,V)$ such that
$i_\beta^\ast(\delta_\beta(e))(h)\ne0$ for all diagrams $\beta$ and edges $e$ of
$\beta$. This exists since we only exclude the union of finitely many
codimension 1 subspaces. By construction, the perturbation represented by $h$
gives every edge of every diagram in the algorithm a color with a nontrivial
$\varepsilon$-component. In particular, all colors are typical.
\end{proof}

Let $\{T_i\}$ be the set of reduced trees. So far, we proved the following for
1-particle irreducible webs. We complete the argument for 1-particle reducible
webs.

\begin{proposition}\label{prop-reduce}
For any connected planar $\gsl$-web $\alpha$ in $(B^3, V_{1,n})$, there is
a perturbation $\alpha(\varepsilon)$ that reduces to a linear combination of
perturbed reduced trees $T_i(\varepsilon)$ in the perturbed skein module
$\skein_\varepsilon(B^3,V_\varepsilon)$
\begin{equation}
\alpha(\varepsilon)=\sum_i c_i(\varepsilon)T_i(\varepsilon).
\end{equation}
\end{proposition}

\begin{proof}
Suppose $\alpha$ is 1-particle reducible. Recall that $\alpha$ must contain the
piece in Figure~\ref{fig-1PR}. We can find a minimal piece of this form with
respect to inclusion, where $\beta$ in the figure is necessarily 1-particle
irreducible. Since $k$ is typical, $a\ne b$ regardless of the sign of the
vertex.

Now apply the reduction algorithm to $\beta\in\skein(B^3,V_{1,1})$ after a
perturbation. The only tree in $(B^3,V_{1,1})$ is a single edge, which would
require $a=b$, contradicting typical colors as above. This means in the
reduction, one of the bubble move \eqref{eq-bubble} must be of the same sign,
which implies that $\beta(\varepsilon)=0$ in the perturbed skein module.
Extending the perturbation to $\alpha$, this proves the proposition by taking
$c_i(\varepsilon)=0$.
\end{proof}

\subsection{Grid paths and trees}
\label{sec-path-tree}

In this section, we ignore edge colors but keep vertex signs to focus on the
combinatorics.

So far, we have assumed the boundary points $V_{1,n}$ consist of one point at
the bottom and $n$ at the top. For the purpose of induction, it is convenient to
consider more general trees with $m$ points at the bottom and $n$ points at the
top. 

We only consider half-polygons with a top boundary edge. We say $T$ is a top
tree if every vertex belongs to a half-polygon. Trivially, a tree with one
bottom point is a top tree. The sign list of a top tree $T$ is a sequence
$S_1,\dotsc,S_{n-1}$ with $\sum\abs{S_i}=n+m-2$, obtained by the same process as
before. A top tree is reduced if every entry is reduced. Let
$\trees^{(n)}_{k,\ell}$ denote the set of reduced top trees (without edge
colors) with $k+1$ points at the bottom, $n$ points at the top, and $k+\ell$
vertices colored $-$.

Recall the set of grid paths $\paths^{(n)}_{k,\ell}$ from Section~\ref{sub.mult}
and the moves $M_s$ from \eqref{eq-moves}. We define a map
\begin{equation}
f^{(n)}_{k,\ell}:\paths^{(n)}_{k,\ell}\to\trees^{(n)}_{k,\ell}
\end{equation}
by mapping the moves $M_1,M_2,M_3,M_4$ to reduced entries $+-,+,-,\emptyset$ of
sign lists.

We can also construct the trees directly using the following induction. The base
case $P=[\,]$ maps to the identity %
$\begin{linkdiag}
\fill[gray!20] (0,0)rectangle(1,1);
\draw[thick] (0.5,1) -- (0.5,0);
\end{linkdiag}$. %
Now suppose $P\in\paths^{(n)}_{k,\ell}$ and $P'=[P,M_{s_n}]$. Let $T$ be the
reduced top tree corresponding to $P$, and construct the tree for $P'$ by the
following.
\begin{equation}\label{eq-path-tree}
M_1\to
\begin{linkdiag}[0.9]
\fill[gray!20] (0,0) rectangle (2.1,1.8);
\draw[fill=white] (0.1,0.8) rectangle (1.3,1.4) (0.7,1.1)node{$T$};
\webnode{1.2,0.4}{A}{+};
\webnode{1.8,0.4}{B}{-};
\begin{scope}[thick]
\draw (0.2,0) -- +(0,0.8) (0.8,0) -- +(0,0.8) (0.5,0.4)node{...};
\draw (0.2,1.4) -- (0.2,1.8) (1.2,1.4) -- (1.2,1.8) (0.7,1.6) node{...};
\draw (1.2,0) -- (A) -- (1.2,0.8) (A) -- (B) (1.8,0) -- (B) -- (1.8,1.8);
\end{scope}
\end{linkdiag},\quad
M_2\to
\begin{linkdiag}[0.9]
\fill[gray!20] (0,0) rectangle (2,1.8);
\draw[fill=white] (0.1,0.8) rectangle (1.3,1.4) (0.7,1.1)node{$T$};
\webnode{1.7,0.4}{B}{+};
\begin{scope}[thick]
\draw (0.2,0) -- +(0,0.8) (0.8,0) -- +(0,0.8) (0.5,0.4)node{...};
\draw (0.2,1.4) -- (0.2,1.8) (1.2,1.4) -- (1.2,1.8) (0.7,1.6) node{...};
\draw (1.2,0.8) ..controls+(0,-0.4).. (B) (1.7,0) -- (B) -- (1.7,1.8);
\end{scope}
\end{linkdiag},\quad
M_3\to
\begin{linkdiag}[0.9]
\fill[gray!20] (0,0) rectangle (2,1.8);
\draw[fill=white] (0.1,0.8) rectangle (1.3,1.4) (0.7,1.1)node{$T$};
\webnode{1.7,0.4}{B}{-};
\begin{scope}[thick]
\draw (0.2,0) -- +(0,0.8) (0.8,0) -- +(0,0.8) (0.5,0.4)node{...};
\draw (0.2,1.4) -- (0.2,1.8) (1.2,1.4) -- (1.2,1.8) (0.7,1.6) node{...};
\draw (1.2,0.8) ..controls+(0,-0.4).. (B) (1.7,0) -- (B) -- (1.7,1.8);
\end{scope}
\end{linkdiag},\quad
M_4\to
\begin{linkdiag}[0.9]
\fill[gray!20] (0,0) rectangle (2,1.8);
\draw[fill=white] (0.1,0.8) rectangle (1.3,1.4) (0.7,1.1)node{$T$};
\begin{scope}[thick]
\draw (0.2,0) -- +(0,0.8) (0.8,0) -- +(0,0.8) (0.5,0.4)node{...};
\draw (0.2,1.4) -- (0.2,1.8) (1.2,1.4) -- (1.2,1.8) (0.7,1.6) node{...};
\draw (1.2,0.8)..controls+(0,-0.4) and +(0,-0.4)..(1.7,0.8) -- (1.7,1.8);
\end{scope}
\end{linkdiag}.
\end{equation}
Note the grid path condition is used to ensure there is at least one strand at
the bottom for the construction. It is clear that the new tree is also a reduced
top tree.

Conversely, every reduced top tree looks like one of the options in
\eqref{eq-path-tree}. A similar induction constructs the inverse of
$f^{(n)}_{k,\ell}$. This proves $f^{(n)}_{k,\ell}$ is a bijection. In
particular, we have the following.

\begin{lemma}\label{lemma-tree-count}
The number of reduced trees in $(B^3, V_{1,n})$ is
$m^{(n)}_{0,\ell}$, the expected multiplicity from \eqref{otimesn}.
\end{lemma}

\begin{example}
Let $P=[M_2,M_1,M_3,M_4]$, which is a grid path ending at $(0,2)$. The path and
its corresponding tree in $(B^3, V_{1,5})$ is shown in
Figure~\ref{fig-path-tree}.

\begin{figure}
\centering
\begin{tikzpicture}[baseline=1cm]
\draw[gray] (0,0) grid (2,2);
\begin{scope}[thick]
\draw[-o-={0.2}{>}] (0,0) to[out=90,in=135] node[pos=0.4,left]{$M_2$} (0.5,0.5) to[out=-45,in=0] (0,0);
\draw[-o-={0.6}{>}] (0,0) -- node[below]{$M_1$} (1,0);
\draw[-o-={0.6}{>}] (1,0) -- node[right]{$M_3$} (1,1);
\draw[-o-={0.6}{>}] (1,1) -- node[above right,inner sep=2pt,pos=0.4,fill=white]{$M_4$} (0,2);
\end{scope}
\end{tikzpicture}
\qquad
\begin{tikzpicture}[baseline=2.5cm]
\fill[gray!20] (-0.5,0.8) rectangle (2.5,4);
\webnode{0.5,3.6}{M1}{+};
\webnode{0.5,2.8}{M2L}{+};
\webnode{1,2.8}{M2R}{-};
\webnode{1.5,2}{M3}{-};
\begin{scope}[thick]
\draw (0,4)..controls+(0,-0.4)..(M1);
\draw (0.5,0.8) -- (M2L) -- (M1) -- (0.5,4) (M2L) -- (M2R);
\draw (1,4) -- (M2R) -- (1,2.4)..controls+(0,-0.4)..(M3);
\draw (1.5,4) -- (M3) -- (1.5,1.6) arc[radius=0.25,start angle=-180,end angle=0] -- (2,4);
\end{scope}
\draw[red] (0.8,4) |- (-0.5,3.2) (1.3,4) |- (-0.5,2.4) (1.8,4) |- (-0.5,1.6);
\end{tikzpicture}
\caption{The grid path $P=[M_2,M_1,M_3,M_4]$ and the corresponding reduced tree.}
\label{fig-path-tree}
\end{figure}
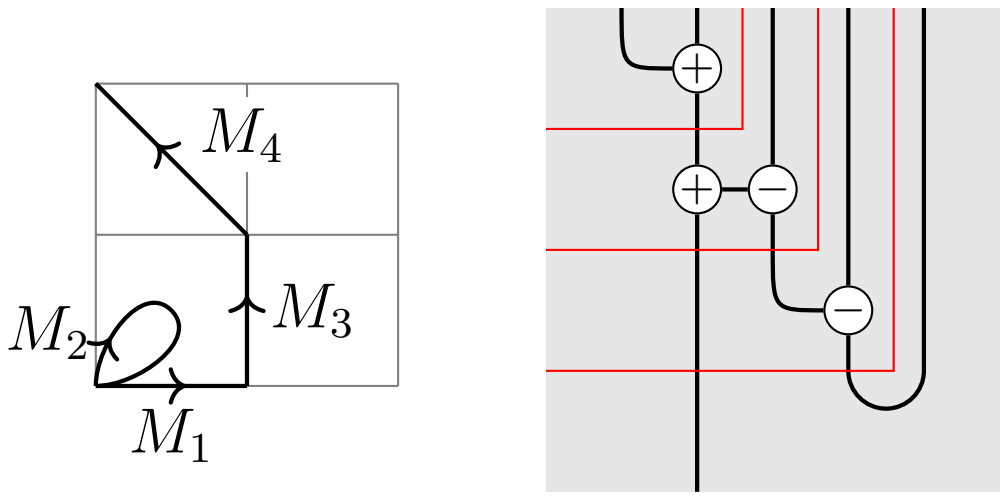
\end{example}

\subsection{Highest terms of trees}
\label{sec-sk-highest}

As in the last section, we also generalize the map $F_{1,n}$ from \eqref{eq-F1n}
to $F_{m,n}$ with $m$ points at the bottom and $n$ points at the top.

\begin{lemma}\label{lemma-tree-high}
Let $P=[M_{s_1},\dotsc,M_{s_{n-1}}]\in\paths^{(n)}_{k,\ell}$ be a grid path and
$T=f^{(n)}_{k,\ell}(P)$ be the corresponding reduced top tree. Suppose the
boundary colors are such that the unique coloring of $T$ is typical. Then in
the alternative height order,
\begin{equation}
F_{k+1,n}(T)(v_1\otimes\dotsm\otimes v_1)=c_P v_1\otimes v_{s_1}\otimes\dotsm\otimes v_{s_{n-1}}+\text{lower terms}
\end{equation}
with a nonvanishing leading coefficient $c_P\ne0$.
\end{lemma}

\begin{proof}
The base case $P=[\,]$ is trivial. For the inductive step, suppose the lemma is
true for paths of length $n-1$. Let
$P'=[M_{s_1},\dotsc,M_{s_{n-1}},M_{s_n}]\in\paths^{(n+1)}_{k',\ell'}$ and
$T'=f^{(n+1)}_{k',\ell'}(P')$ whose top and bottom colors are
$a_1,\dotsc,a_{n+1}$ and $b'_1,\dotsc,b'_{k'+1}$. Then
$P=[M_{s_1},\dotsc,M_{s_{n-1}}]\in\paths^{(n)}_{k,\ell}$ satisfy the lemma, and
write $T=f^{(n)}_{k,\ell}(P)$. By construction \eqref{eq-path-tree}, $T$ appears
as the top left part of $T'$, so the top and bottom colors of $T$ are
$a_1,\dotsc,a_n$ and $b_1,\dotsc,b_{k+1}$ where $b_i=b'_i$ for $i=1,\dotsc,k$.
Then
\begin{equation}
\begin{split}
F_{k'+1,n+1}(T')&:\bigotimes_{i=1}^{k'+1} V(0,b'_i)
=\Big(\bigotimes_{i=1}^k V(0,b_i)\Big)\otimes\Big(\bigotimes_{i=k+1}^{k'+1}V(0,b'_i)\Big)\\
&\quad\xrightarrow{\id\otimes g_{s_n}}
\Big(\bigotimes_{i=1}^k V(0,b_i)\Big)\otimes\big(V(0,b_{k+1})\otimes V(0,a_{n+1})\big)\\
&\quad=\Big(\bigotimes_{i=1}^{k+1} V(0,b_i)\Big)\otimes V(0,a_{n+1})
\xrightarrow{F_{k+1,n}(T)\otimes\id}\bigotimes_{i=1}^n V(0,a_i).
\end{split}
\end{equation}
Here, $g_{s_n}$ is the map corresponding to the added part in the lower right
corner of \eqref{eq-path-tree}. By definition,
\begin{equation}
g_{s_n}=\begin{cases}
(\id\otimes d_e\otimes\id)(\gamma^{b_{k+1},e^\ast}_+\otimes\gamma^{e,a_{n+1}}_-),& s_n=1,\\
\gamma^{b_{k+1},a_{n+1}}_+,& s_n=2,\\
\gamma^{b_{k+1},a_{n+1}}_-,& s_n=3,\\
b_{b_{k+1}},& s_n=4,
\end{cases}
\end{equation}
where $e=b_{k'+1}-a_{n+1}-1$, and the basic maps are defined in
Section~\ref{sec-basic-maps}. A straightforward calculation shows that
$g_{s_n}(\bigotimes_{i=k+1}^{k'+1}v_1)=c_{s_n}v_1\otimes v_{s_n}+(\text{lower terms})$
with $c_{s_n}\ne0$ if all colors are typical. The $s_n=1$ case is given in
\eqref{eq-Hpm-highest}, while the others directly come from the formulas in
Section~\ref{sec-basic-maps}. Then by induction,
\begin{equation}
\begin{split}
F_{k'+1,n+1}(T')(v_1\otimes\dotsm\otimes v_1)
&=c_{s_n}(F_{k+1,n}(T)\otimes\id)\Big(\big(\bigotimes_{i=1}^k v_1\big)\otimes v_1\otimes v_{s_n}+\text{lower terms}\Big)\\
&=c_{s_n}c_P v_1\otimes v_{s_1}\otimes\dotsm\otimes v_{s_{n-1}}\otimes v_{s_n}+\text{lower terms}.
\end{split}
\end{equation}
In the second step, the image of the lower terms are lower in the
lexicographical order by considering the weights.
\end{proof}

\begin{remark}
As mentioned in Remark~\ref{rem-basis}, the alternative height order is chosen
for the proof of this lemma. We elaborate on the reasons here.

The inductive construction of $f^{(n)}_{k,\ell}$ is parallel to the definition
\eqref{eq-path-incl} of $i_P$, so the proof is parallel to that of
Lemma~\ref{lemma-highest}, with $\bigotimes_{i=1}^{k+1} V(0,b_i)$ playing the
role of $V_P=V(k,\sum_{i=1}^{k+1}b_i)$, and $v_1\otimes\dotsm\otimes v_1$ is the
highest weight vector of $V_P$ embedded in the tensor product.

However, in the traditional height order, the highest component is
$V(0,\sum_{i=1}^{k+1}b_i)$ instead, and many more components are higher than
$V_P$ in general. It is difficult to remove these higher components for the
induction to work.
\end{remark}

\subsection{Proof of Theorem~\ref{thm.1}}
\label{sec-proof}

Now we go back to $F_{1,n}$ with only one point on the bottom. Recall the
definition of a generic list from Section~\ref{sub.mult}. We start by assuming
the top colors $a_1,\dotsc,a_n$ are generic. The following result is a simple
combinatorial calculation.

\begin{lemma}\label{lemma-typ-tree}
Suppose $T$ is a trivalent tree in $(B^3,V_{1,n})$ with vertex signs. If the top
colors $a_1,\dotsc,a_n$ are generic, then the unique coloring of $T$ is typical.
\end{lemma}

\begin{corollary}\label{cor-indep}
If the top colors $a_1,\dotsc,a_n$ are generic, then the reduced trees in
$\skein(B^3,V_{1,n})$ are linearly independent.
\end{corollary}

\begin{proof}
The conditions of Lemma~\ref{lemma-tree-high} are satisfied by
Lemma~\ref{lemma-typ-tree}. Comparing with Lemma~\ref{lemma-highest}, the images
of reduced trees by $F_{1,n}$ is a basis, so the trees are linearly
independent.
\end{proof}

Let $\{T_P=f^{(n)}_{0,\ell}(P)\}$ denote the set of reduced trees. We want to
show that they span the skein module $\skein(B^3,V_{1,n})$. Recall that we can
always resolve crossings, so we only need to consider planar webs. Given a
connected planar $\gsl$-web $\alpha$, by Proposition~\ref{prop-reduce}, there
exists a perturbation $\alpha(\varepsilon)$ which can be expanded
\begin{equation}\label{eq-perturb-span}
\alpha(\varepsilon)=\sum_P c_P(\varepsilon)T_P(\varepsilon)
\end{equation}
in the perturbed skein module $\skein_\varepsilon(B^3,V_\varepsilon)$. Note
$c_P(\varepsilon)\in R_\varepsilon$. If some $c_P(\varepsilon)$ is not finite at
$\varepsilon=0$, we can multiply both sides by the appropriate power of
$\CC{\varepsilon}$ to cancel the poles exactly. Then we can set $\varepsilon=0$
to conclude that the trees $\{T_P\}$ are linearly dependent.

If the top colors $a_1,\dotsc,a_n$ are generic, then by
Corollary~\ref{cor-indep}, $\{T_P\}$ is in fact linearly independent, which
means that all coefficients $c_P(\varepsilon)$ are finite by contradiction. Then
setting $\varepsilon=0$ in \eqref{eq-perturb-span}, every connected $\gsl$-web is in the $\groundR$-span of $\{T_P\}$.

Finally, to deal with disconnected webs, apply the same argument to each connected component. By a nesting argument, there is a tree whose boundary points is a consecutive subset of the top colors. Then the conservation in this component shows that the top colors are not generic, a contradiction.

This proves $\{T_P\}$ is a basis of the skein module $\skein(B^3,V_{1,n})$ over
$\groundR$ for generic colors, and the rank is the expected number by
Lemma~\ref{lemma-tree-count}.

Now suppose the colors $\vec{a}=(a_1,\dotsc,a_n)$ are not generic.
Motivated by the tensor rule of $V(0,a_i)$, we prove this case by establishing
an isomorphism $\vec{a}\cong\vec{a}'$ of objects of $\skeinballcat$ where
$\vec{a}'=(a'_1,\dotsc,a'_n)$ is generic.

We start with the simplest nontrivial case.

\begin{lemma}\label{lemma-sk-iso} %
Suppose $a,b,c,d$ are typical colors with $a+b=c+d$ and $a-c$ typical. Then the
skein
\begin{equation}
H_{a,b}^{c,d}=\CC{c-a}\Hskein{+}{-}
\end{equation}
is an isomorphism $(a,b)\xrightarrow{\cong}(c,d)$ in $\skein(B^3)$.
Moreover,
\begin{equation}\label{eq-Hpm-highest}
F(H_{a,b}^{c,d})(v_1\otimes v_1)=v_1\otimes v_1.
\end{equation}
\end{lemma}

\begin{proof}
The inverse is $H_{c,d}^{a,b}$ by square breaking \eqref{eq-square}.
\eqref{eq-Hpm-highest} is a straightforward calculation.
\end{proof}

For the general case, we can build an isomorphism out of a chain of $H$'s. If we
choose $\vec{a}'$ carefully, the typical condition of the lemma is always
satisfied. This finishes the proof of Theorem~\ref{thm.1}.\hfill\qed

\begin{remark}\label{rem-other-bases} %
Note that the proof gives a basis in the non-generic case by applying a fixed
isomorphism to a generic tree basis. This is used in Example~\ref{exa-basis2} to
obtain Figure~\ref{fig-basis-end}. However, in practice, we can find a set of
webs with the same leading term property as Lemma~\ref{lemma-tree-high}, then
conclude a basis by dimension count. Of course, this only works a posteriori
knowing the rank of the skein module in the non-generic case, since the spanning
property is not clear.

A good example is $\End_\skeinballcat(a_1,\dotsc,a_n)=\skein(B^3,V_{n,n})$ with
$n\ge2$. (The case when all $a_i=a$ is the focus of Section~\ref{sub.Bn}.) By
duality, this is $\skein(B^3,V_{1,2n-1})$ with top colors
$(a_1,\dotsc,a_n,a_n^\ast,\dotsc,a_2^\ast)$, which is necessarily non-generic.
For the first grid path
$P=[M_1,\dotsc,M_1,M_4,\dotsc,M_4]\in\paths^{(2n-1)}_{0,n-1}$, the expected
leading term from Lemma~\ref{lemma-tree-high} is
\begin{equation}
F_{1,2n-1}(T_P)(v_1)=c_Pv_1\otimes\dotsm\otimes v_1\otimes v_4\otimes\dotsm\otimes v_4+\text{lower terms}.
\end{equation}
On the other hand, it is clear by \eqref{eq-ba} that the identity in
$\End_\skeinballcat(a_1,\dotsc,a_n)$ has the same leading term after duality.
Therefore, we can replace the first basis element with the identity diagram.
\end{remark}

%%%%%%%%%%%%%%%%%%%%%%%%%%%%%%%%%%%%%%%%%%%%%%%%%%%%%%%%%%%%%%%%%%%%%%%%%%%%
%%%%%%%%%%%%%%%%%%%%%%%%%%%%%%%%%%%%%%%%%%%%%%%%%%%%%%%%%%%%%%%%%%%%%%%%%%%%

\section{The colored Links--Gould polynomials of links}

\subsection{Definition}
\label{sub.LGdef}

In this section we recall the definition of the colored
Links--Gould polynomial of a link. The following is essentially 
proved by \cite{GP:multi}.

\begin{proposition}\label{prop-S3}
$\skein(S^3)$ is a rank 1 free $\groundR$-module generated by the
unknot $U_a$ with any typical color $a\in\colorL\setminus\{0,-1\}$. Moreover,
$\CC{a}\CC{a+1}U_a$ is independent of the color $a$.
\end{proposition}

Suppose $L\subset S^3$ is an oriented link with $r$ ordered components. Let
$\tilde{L}$ be $L$ with $0$-framing for simplicity and a color $a_i$ on the
$i$-th component. Define the $\gsl$-bracket $\bracket{\cdot}_a$ of $\tilde{L}$
with respect to a typical reference color $a$ by
\begin{equation}
\tilde{L}=\bracket{\tilde{L}}_a\cdot U_a\in\skein(S^3).
\end{equation}
Let $\mu_{i,j}$ be the linking matrix of $\tilde{L}$, with the convention that
the diagonal elements are the framings, which are $0$. Then there exists a
normalized invariant $\GG_L(q;t_1,\dotsc,t_r)$, independent of the reference
color $a$, such that for any choice of colors $a_1,\dotsc,a_r$ (not necessarily
the basis colors of $\colorL$),
\begin{equation}
\label{eq-ML}
\GG_L(q;q^{a_1},\dotsc,q^{a_r})=\frac{q^w}{\CC{a}\CC{a+1}}\bracket{\tilde{L}}_a,\qquad
w=\sum_{i>j}\mu_{i,j}(2a_i+1)(2a_j+1),
\end{equation}
where $w$ is the writhe correction similar to the Jones polynomial via Kauffman
brackets. Note in general, there is another term in $w$ involving framings.

Finally, define the multivariable colored Links--Gould polynomial of $L$ as
\begin{equation}\label{eq-LG}
\LG_L(q;t_1,\dotsc,t_r)=\begin{cases}
(t_1-t_1^{-1})(qt_1-q^{-1}t_1^{-1})\GG_L(q;t_1),& r=1,\\
\GG_L(q;t_1,\dotsc,t_r),& r\ge2
\end{cases}
\end{equation}
with $t_i=q^{a_i}$. 
Note that $\LG_L(q;t_1,\dotsc,t_r) 
\in \BZ[q^{\pm 1},t_1^{\pm 1},\dotsc,t_r^{\pm 1}]$, and that
the usual Links--Gould polynomial~\cite{Links-Gould} can be recovered by setting all $t_i$ to be equal.

\begin{example}
Since many conventions and normalizations are used in the literature, we give a
few examples for clarity.
\begin{equation}
\begin{aligned}
\LG_{T(2,1)}(q;t_1)&=1,\\
\LG_{T(2,2)}(q;t_1,t_2)&=1,\\
\LG_{T(2,3)}(q;t_1)&=t_1^4-(1+q^{-2})t_1^2+(1+2q^{-2})-(q^{-4}+q^{-2})t_1^{-2}+q^{-4}t_1^{-4},\\
\LG_{T(2,4)}(q;t_1,t_2)&=qt_1^2t_2^2+(q+q^{-1})+q^{-3}t_1^{-2}t_2^{-2}.
\end{aligned}
\end{equation}
Here, $T(p,q)$ denotes the $(p,q)$ torus link. Specifically, $T(2,1)$ is the
unknot, $T(2,2)$ is the Hopf link with linking number $+1$, $T(2,3)$ is the
right-handed trefoil knot, and $T(2,4)$ is the two-component link with crossing
number $4$ and linking number $+2$.
\end{example}

We now consider the higher Links--Gould polynomials
$\LG^{(n_1,\dotsc,n_r)}_L(q;t_1,\dotsc,t_r)$, which depend on a list of integers
$n_i\ge1$, where the basic Links--Gould polynomial $\LG_L$ is the specialization
$n_i=1$. Recall the ideal $\mathcal{I}$ of projective objects in $\Repsl$
discussed in Remark~\ref{rem.proj}. Using the admissible skein module
$\skein_\mathcal{I}(S^3)$, which is isomorphic to
$\skein(S^3)\otimes_{\groundR}\groundK$ by Theorem~\ref{thm.2}, we can extend
the bracket $\bracket{\cdot}_{m,b}$ to links colored by any projective objects,
which include all typical $V(n-1,a)$ by Lemma~\ref{lem-typ}, with respect to the unknot $U_{m,b}$ colored by $V(m-1,b)$. Define
\begin{equation}\label{eq-higher-ML}
\begin{gathered}
\GG^{(n_1,\dotsc,n_r)}_L(q;q^{a_1},\dotsc,q^{a_r})
=\frac{q^w(-1)^{m-1}\CC{m}}{\CC{1}\CC{b}\CC{b+m}}\bracket{\tilde{L}}_{m,b},\\
w=\sum_{i>j}\mu_{i,j}(2a_i+n_i)(2a_j+n_j),
\end{gathered}
\end{equation}
where the $i$-th component of $\tilde{L}$ is colored by $V(n_i-1,a_i)$. Clearly,
\eqref{eq-higher-ML} specializes to $\GG_L$ for $n_i=1$. As before,
$\GG^{(n_1,\dotsc,n_r)}_L$ is the normalized invariant, and the Links--Gould
invariants are given by
\begin{equation}\label{eq-higer-LG}
\LG^{(n_1,\dotsc,n_r)}_L(q;t_1,\dotsc,t_r)=\begin{cases}
\GG^{(n_1)}_L(q;t_1)/B_{n_1}(q;t_1),& r=1,\\
\GG^{(n_1,\dotsc,n_r)}_L(q;t_1,\dotsc,t_r),& r\ge2,
\end{cases}
\end{equation}
where
\begin{equation}
B_n(q;t)=
(-1)^{n-1}\frac{q^n-q^{-n}}{(q-q^{-1})(t-t^{-1})(q^nt-q^{-n}t^{-1})}
\end{equation}
is the normalized invariant of the unknot.

\subsection{Proof of the cabling theorems}

In this section we prove the cabling theorems for the colored Links--Gould
polynomial of links. We fix a reference color for the brackets and omit it from
notations.

The first proof below is not phrased using the functoriality
(Lemma~\ref{lem-ann-embed}), although it is clear how it applies.

\begin{proof}(of Theorem~\ref{thm-ML-recursion})
Choose $\tilde{L}$ as in the definition of $\GG_L$. The bracket
$\bracket{\cdot}$ satisfies a simple decomposition formula following
\eqref{otimesn} (see e.g.~\cite[Thm.~1.2]{Ga:colored})
\begin{equation}
\bracket{\tilde{L}(p+1,0)}
=\sum_{k,\ell}m^{(p+1)}_{k,\ell}\bracket{\tilde{L}^{(k+1)}_{b+\ell}}.
\end{equation}
Here, $\tilde{L}(p+1,0)$ is $\tilde{L}$ with the last component replaced by
$p+1$ parallel copies colored by $V(0,a_r),\dotsc,V(0,a_{r+p})$, and
$\tilde{L}^{(k+1)}_{b+\ell}$ is $\tilde{L}$ except the last component is colored
by $V(k,b+\ell)$ where $b=\sum_{i=r}^{r+p}a_r$. Then by a careful calculation of
the writhe corrections, we get \eqref{eq-ML-recursion}.
\end{proof}

\begin{example}\label{exa-Hopf} %
Suppose $L=T(2,2)$ is the Hopf link with linking number $+1$. Using
\cite[Lemma~3.15]{Ga:colored}, we can calculate
\begin{equation}
\LG^{(n_1,n_2,\dotsc,n_{p+2})}_{L(p+1,0)}(q;t_1,\dotsc)
=(-1)^{\sum_{i=1}^{p+2}(n_i-1)}\Big(\frac{\CC{1}}{\CC{n_1}}d_{n_1}(q;t_1)\Big)^p
\prod_{i=2}^{p+2}\frac{\CC{n_1n_i}}{\CC{1}},
\end{equation}
where
\begin{equation}
d_n(q;t_1)=(t_1-t_1^{-1})(q^nt_1-q^{-n}t_1^{-1}).
\end{equation}
In particular,
\begin{align}
\LG^{(n_1,n_2)}_L(q;t_1,t_2)
&=(-1)^{n_1+n_2}\frac{\CC{n_1n_2}}{\CC{1}},\\
\LG^{(n_1,1,\dotsc,1)}_{L(p+1,0)}(q;t_1,\dotsc)
&=(-1)^{n_1-1}\frac{\CC{n_1}}{\CC{1}}d_{n_1}(q;t_1)^p.
\end{align}
It is straightforward to verify \eqref{eq-ML-recursion} using the generating
function \eqref{eq-gen-fn}.
\end{example}

\begin{proof}(of Theorem~\ref{thm-ML-cable}) %
By functoriality, we just need to prove the corresponding statement in the
annulus $A$. Let $c$ denote the core of the annulus, $c_b$ denote $c$ colored by
$b\in\colorL$, and $c(n,1)$ denote the $(n,1)$-cable of $c$ in the thickened
annulus $A\times[0,1]$ colored by $a$. Further, define two families of webs
$\alpha^{(n)}_{r,\ell},\beta^{(n)}_{r,\ell}\in\skein(A),r=1,\dotsc,n-1$ by the
closure of the following tangles with a preemptive writhe correction
$w^{(n)}_{r,\ell}=q^{(n-r-1)(2a^2+2a)}\CC{ra+\ell}$.
\begin{equation}
\alpha^{(n)}_{r,\ell}=-w^{(n)}_{r,\ell}\,
\begin{tikzpicture}[baseline=(ref.base)]
\fill[gray!20] (-0.6,0) rectangle (2,2);
\webnode{0,1.5}{A}{-};
\webnode{0,0.5}{B}{+};
\begin{scope}[thick]
\draw (A) to[out=-135,in=135] (B) (0,0) -- (B) (A) -- (0,2);
\begin{knot}
\strand (1.6,2) ..controls (1.6,1.6).. (B);
\strand (A) ..controls +(0.4,-0.4).. (0.4,0)
  (0.4,2) ..controls (0.4,1.5) and (0.8,1.5).. (0.8,0)
  (1.2,2) ..controls (1.2,1.5) and (1.6,1.5).. (1.6,0);
\end{knot}
\path[-o-={0.8}{>}] (B) -- (0,0);
\path[-o-={0.5}{>}] (1.6,2) ..controls (1.6,1.6).. (B);
\path[-o-={0.8}{>}] (A) ..controls +(0.4,-0.4).. (0.4,0);
\path[-o-={0.8}{>}] (0.4,2) ..controls (0.4,1.5) and (0.8,1.5).. (0.8,0);
\path[-o-={0.8}{>}] (1.2,2) ..controls (1.2,1.5) and (1.6,1.5).. (1.6,0);
\end{scope}
\path (1.2,0.3)node{...};
\draw[decoration=brace,decorate] (1.6,-0.1) --node[below]{$n-r-1$} 
(0.4,-0.1);
\node (ref) at (0,1) {\phantom{$-$}};
\end{tikzpicture}\,,
\qquad
\beta^{(n)}_{r,\ell}=w^{(n)}_{r,\ell}\,
\begin{tikzpicture}[baseline=(ref.base)]
\fill[gray!20] (-0.6,0) rectangle (2,2);
\webnode{0,1.5}{A}{+};
\webnode{0,0.5}{B}{-};
\begin{scope}[thick]
\draw (A) to[out=-135,in=135] (B) (0,0) -- (B) (A) -- (0,2);
\begin{knot}
\strand (1.6,2) ..controls (1.6,1.6).. (B);
\strand (A) ..controls +(0.4,-0.4).. (0.4,0)
  (0.4,2) ..controls (0.4,1.5) and (0.8,1.5).. (0.8,0)
  (1.2,2) ..controls (1.2,1.5) and (1.6,1.5).. (1.6,0);
\end{knot}
\path[-o-={0.8}{>}] (B) -- (0,0);
\path[-o-={0.5}{>}] (1.6,2) ..controls (1.6,1.6).. (B);
\path[-o-={0.8}{>}] (A) ..controls +(0.4,-0.4).. (0.4,0);
\path[-o-={0.8}{>}] (0.4,2) ..controls (0.4,1.5) and (0.8,1.5).. (0.8,0);
\path[-o-={0.8}{>}] (1.2,2) ..controls (1.2,1.5) and (1.6,1.5).. (1.6,0);
\end{scope}
\path (1.2,0.3)node{...};
\draw[decoration=brace,decorate] (1.6,-0.1) --node[below]{$n-r-1$} 
(0.4,-0.1);
\node (ref) at (0,1) {\phantom{$-$}};
\end{tikzpicture}\,.
\end{equation}
Here, all crossings have equal color $a$, and the outward labeled edge of the
bottom vertex is colored $(r+1)a+\ell$ in both $\alpha$ and $\beta$. When
$r=n-1$, both diagrams consist of a single bubble. When $r=0$, the definitions
are extended by $\alpha^{(n)}_{0,0}=\beta^{(n)}_{0,0}=q^{(n-1)(2a^2+2a)}c(n,1)$
based on \eqref{eq-0break}. 

By resolving the leftmost crossing using \eqref{eq-equal-cross}, the first 
term
is a disconnected diagram with a bubble that can be reduced to a product
$c_{(r+1)a+l}\alpha^{(n-r-1)}_{0,0}$. The other two terms contain squares, 
which can be
reduced by the algorithm in Section~\ref{sec-reduce}, replacing it locally 
with
$\Iskein{\pm}{\mp}$. We can then pull the bottom vertex around the annulus,
creating new $\alpha,\beta$ diagrams with $r$ increased by one. After a 
careful
calculation, the result is
\begin{equation}
\begin{aligned}
\alpha^{(n)}_{r,\ell}&=
-q^{-ra-\ell}c_{(r+1)a+\ell}\alpha^{(n-r-1)}_{0,0}+q^a\alpha^{(n)}_{r+1,\ell}
+q^{-(r+1)a-\ell-1}\frac{\CC{a+1}}{\CC{(r+1)a+\ell+1}}\beta^{(n)}_{r+1,\ell+1},\\
\beta^{(n)}_{r,\ell}&=
-q^{ra+\ell}c_{(r+1)a+\ell}\alpha^{(n-r-1)}_{0,0}
+q^{(r+1)a+\ell}\frac{\CC{a}}{\CC{ra+\ell}}\alpha^{(n)}_{r+1,\ell}
+q^{-a-1}\beta^{(n)}_{r+1,\ell+1},
\end{aligned}
\end{equation}
with initial conditions
\begin{equation}
\alpha^{(n)}_{n-1,\ell}=\frac{\CC{na+\ell}}{\CC{a}}c_{na+\ell},\qquad
\beta^{(n)}_{n-1,\ell}=\frac{\CC{na+\ell+1}}{\CC{a+1}}c_{na+\ell}.
\end{equation}
The solutions are given by
\begin{multline}\label{eq-alpha-sol}
(-1)^{n-1}\alpha^{(n)}_{r,\ell}
=(-1)^rq^{-ra-l} c^{(n-r)}_{na+\ell}\\
+\sum_{i=r+1}^{n-1}(-1)^i\Big(q^{(2i-3r)a-l}c^{(n-i)}_{na+\ell}
+\frac{q^{-2i(a+1)+2(r-\ell)}\CC{a+1}}{{\CC{(r+1)a+\ell+1}}}c^{(n-i)}_{na+i-r+\ell}\Big)\\
+\CC{ra+\ell}\Big((-1)^{n-1}\frac{q^{(2n-2r-1)a}}{\CC{a}}c^{(1)}_{na+\ell}
+\sum_{i=r+2}^{n-1}(-1)^{i+1}\frac{q^{(2i-3r-3)a-\ell-1}}{\CC{(r+1)a+\ell+1}}c^{(n-i)}_{na+\ell+1}\Big),
\end{multline}
where $c^{(k+1)}_b$ denote $c$ colored by $V(k,b)$. $\beta^{(n)}_{r,\ell}$ is
similar and thus omitted. Specializing to $(r,\ell)=(0,0)$, the last line
vanishes, and we can rearrange the rest to get
\begin{equation}
(-1)^{n-1}q^{(n-1)(2a^2+2a)}c(n,1)
=c^{(n)}_{na}+\sum_{i=1}^{n-1}(-1)^i
\big(q^{2ia} c^{(n-i)}_{na}+q^{-2i(a+1)}c^{(n-i)}_{na+i}\big).
\end{equation}
The theorem follows from a simple calculation of the writhe corrections.
\end{proof}

\subsection{Proof of the specialization at $q=1$}
\label{sec.alex}

In this section we discuss the specialization of the
$\gsl$-skein module at $q=1$ focusing on links in $S^3$. Contrary to what
happens for the skein module of a simple Lie algebra, which is trivial
under this specialization, the $\gsl$-skein module is not.

A key part in the proof that the colored Links--Gould polynomial
of a link specializes when $q=1$ to the square of the multivariable Alexander polynomial is a non-local characterization of the latter,
given by Turaev~\cite[Sec.4]{Tu:alex}. We will follow the
same ideas here to prove Theorem~\ref{thm-alex}, first for the basic colors
$n_i=1$, and then by induction for all colors.

\begin{lemma}
\label{lem.tur}
Suppose $f_L(t_1,\dotsc,t_r)$ is a rational function invariant of oriented links with ordered components that is invariant under simultaneous
permutations of the components and variables and satisfying the cabling
condition
\begin{equation}\label{tor.1}
f_{L(2,1)}(t_1,\dotsc,t_r)+f_{L(2,0)}(t_1,\dotsc,t_r,t_r)
=c\Big(\prod_{j=1}^{r-1}t_j^{\mu_{r,j}},t_r\Big)f_L(t_1,\dotsc,t_{r-1},t_r^2),
\end{equation}
where $\mu_{i,j}$ is the linking matrix of $L$, and
\begin{equation}\label{eq-cable-scalar}
c(T,t_r)=t_r^2T^2+T^2+T^{-2}+t_r^{-2}T^{-2}.
\end{equation}
Then $f$ is uniquely determined by its diagonal specializations
$f_L(t,\dotsc,t)$ on all links $L$.
\end{lemma}

\begin{proof}
Since these conditions are linear, it suffices to show that $f=0$ if the
diagonal specialization is $0$, i.e., $f_L(t,\dotsc,t)=0$ for all links $L$.
Using \eqref{tor.1}, we can show $f_L(t^{2^{k_1}},\dotsc,t^{2^{k_r}})=0$ by an
induction on $(k_1,\dotsc,k_r)\in\nats^r$. The only rational function with this
property is $0$.
\end{proof}

Now we prove that both the Conway potential and the $q=1$ specialization of
normalized invariants $\GG_L$ satisfy the cabling condition in
Lemma~\ref{lem.tur}. For the Conway potential, this follows from the following
results (known collectively in the literature as Torres relations). They were
proved in \cite[Thm.~5.1--5.4]{SW:Alex}. See also \cite[Thm.~1.3.1]{Tu:alex} and
\cite[Cor.~3.4]{Cimasoni}.

\begin{theorem}
\label{thm-alex-cable} 
Let $L\subset S^3$ be an oriented link with $r$ ordered components and
linking matrix $\mu_{i,j}$.
\begin{enumcc}[(1)]
\item Recall $L(p+1,0)$ is the link obtained by adding $p$ parallel copies of
the last component as the $(r+1)$-th to the $(r+p)$-th components. Then
\begin{equation}\label{eq-torres1}
\nabla_{L(p+1,0)}(t_1,\dotsc,t_{r+p}) = (T-T^{-1})^p\nabla_{L}\Big(t_1,\dotsc,t_{r-1},\prod_{i=r}^{r+p}t_i\Big),
\quad T=\prod_{j=1}^{r-1}t_j^{\mu_{r,j}}.
\end{equation}
\item Suppose $p,q$ are coprime integers. Let $L(p,q)$ be the link obtained by
replacing the last component of $L$ by its $(p,q)$-cable. Then
\begin{equation}\label{eq-torres2}
\nabla_{L(p,q)}(t_1,\dotsc,t_r) = \frac{T_q^p-T_q^{-p}}{T_q-T_q^{-1}}\nabla_L(t_1,\dotsc,t_{r-1},t_r^p),
\quad T_q = t_r^qT. 
\end{equation}
\end{enumcc}
\end{theorem}

\begin{proof}(of Theorem~\ref{thm-alex}) 
We will prove Equation~\eqref{eq-higher-alex}.
First we consider the base case $n_i=1$. Clearly, both the Conway potential
function $\nabla_L(t_1,\dotsc,t_r)$ and the $q=1$ specialization
$\GG_L(1;t_1,\dotsc,t_r)$ are permutation invariant. By
Theorems~\ref{thm-ML-cable}, \ref{thm-ML-recursion}, and \ref{thm-alex-cable},
both invariants also have the property \eqref{tor.1}. Finally, the diagonal
specialization of $\GG_L$ at $q=1$ agrees with that of $\nabla_L$ by \cite{Kohli}.
This concludes the base case by Lemma~\ref{lem.tur}.

Since Theorem~\ref{thm-ML-recursion} uniquely determines
$\GG^{(n_1,\dotsc,n_r)}_L$ from $\GG_L$ by induction, we just need to substitute
\eqref{eq-higher-alex} and $q=1$ into \eqref{eq-ML-recursion} and check that the
result is an equality. In other words, we need to verify
\begin{multline}\label{eq-alex-L}
\nabla_{L(n+1,0)}(t_1,\dotsc,t_{r-1},t_r,\dotsc,t_{r+n})^2\\
=\Big(\sum_{k,\ell}m^{(n+1)}_{k,\ell}T^{2n-2k-4\ell}
(-1)^k(k+1)\Big)
\nabla_L\Big(t_1,\dotsc,t_{r-1},\prod_{i=r}^{r+n}t_i\Big)^2.
\end{multline}
Here, the common factor $\prod_{i=1}^{r-1}(-1)^{n_i-1}n_i$ has been removed, and
we have used $Q|_{q=1}=T^2$, where $T=\prod_{i=1}^{r-1}t_i^{\mu_{r,i}}$ is the
same one from \eqref{eq-torres1}. The right-hand side can be evaluated using the
generating function $G(x,y,z)$ of $m^{(n)}_{k,\ell}$ given in~\eqref{eq-gen-fn}.
We have
\begin{equation}
\begin{split}
&\quad\sum_{n=0}^\infty\Big(\sum_{k,\ell}m^{(n+1)}_{k,\ell}T^{-2k-4\ell}(-1)^k(k+1)\Big)x^n\\
&=x^{-1}\frac{\partial(zG(x,y,z))}{\partial z}\Big
|_{\begin{subarray}{l}y=T^{-4}\\z=-T^{-2}\end{subarray}}
=\frac{1}{1-(1-T^{-2})^2x}.
\end{split}
\end{equation}
Expanding the generating function, we get
\begin{equation}
\sum_{k,\ell}m^{(n+1)}_{k,\ell}T^{2n-2k-4\ell}(-1)^k(k+1)
=T^{2n}(1-T^{-2})^{2n}=(T-T^{-1})^{2n}.
\end{equation}
This agrees with the left-hand side of \eqref{eq-alex-L} by
Theorem~\ref{thm-alex-cable}(1).
\end{proof}

%%%%%%%%%%%%%%%%%%%%%%%%%%%%%%%%%%%%%%%%%%%%%%%%%%%%%%%%%%%%%%%%%%%%%%%%%%%%
%%%%%%%%%%%%%%%%%%%%%%%%%%%%%%%%%%%%%%%%%%%%%%%%%%%%%%%%%%%%%%%%%%%%%%%%%%%%

\bibliographystyle{hamsalpha}
\bibliography{biblio}
\end{document}